\documentclass[preprint]{imsart}

\RequirePackage{amsthm,amsmath,amsfonts,amssymb}
\RequirePackage[colorlinks,citecolor=blue,urlcolor=blue]{hyperref}
\RequirePackage[authoryear]{natbib}
\RequirePackage{graphicx}

\usepackage{natbib}
\usepackage{cancel}
\usepackage[dvipsnames]{xcolor}
\usepackage{dsfont}
\usepackage{natbib}
\usepackage{enumitem}
\usepackage{comment}
\usepackage{framed}
\usepackage{multirow}
\usepackage{cleveref}

\setcitestyle{square, comma, numbers,sort&compress, super}

\theoremstyle{plain}
\newtheorem{thm}{Theorem}[section]
\newtheorem*{thm*}{Theorem}

\newtheorem{lem}[thm]{Lemma}

\newtheorem{prop}[thm]{Proposition}

\theoremstyle{definition}

\newtheorem*{remark*}{Remark}

\renewcommand{\epsilon}{\varepsilon}
\renewcommand{\P}{\operatorname{P}}

\newcommand{\E}{\mathsf{E}}
\newcommand{\V}{V}
\newcommand{\Vn}{V_n}

\providecommand{\eps}{\varepsilon}
\renewcommand{\P}{\operatorname{P}}
\newcommand{\Var}{\mathsf{Var}}

\newcommand{\arginf}{\operatornamewithlimits{arginf}}

\newcommand{\argmin}{\operatornamewithlimits{arg\,min}}
\DeclareMathOperator{\divergence}{div}

\newcommand{\R}{\mathbb R}
\newcommand{\X}{\mathcal{X}}

\newcommand{\N}{\mathcal{N}}
\newcommand{\Nb}{\mathcal{N}_{[\;]}}

\newcommand{\G}{\mathcal{G}}
\renewcommand{\H}{\mathcal{H}}
\renewcommand{\L}{{L}}

\newcommand{\e}{\mathcal{E}}
\newcommand{\rv}{\zeta}
\newcommand{\mn}{\|\cdot\|_{L^2(P_n)}}
\newcommand{\m}{\|\cdot\|_{L^2(P)}}
\newcommand{\rmd}{{\, \rm d}}
\newcommand{\rme}{{\rm e}}

\newcommand{\vvvert}[1]{\left\vert\kern-0.25ex\left\vert\kern-0.25ex\left\vert #1 
    \right\vert\kern-0.25ex\right\vert\kern-0.25ex\right\vert}

\begin{document}

\begin{frontmatter}
\title{Empirical Variance Minimization with Applications in Variance Reduction and Optimal Control}
\runtitle{Empirical Variance Minimization}

\begin{aug}
\author[A,B]{\fnms{Denis} \snm{Belomestny}
\ead[label=e1]{denis.belomestny@uni-due.de}},
\author[B]{\fnms{Leonid} \snm{Iosipoi}
\ead[label=e2,mark]{liosipoi@hse.ru}},
\author[B]{\fnms{Quentin} \snm{Paris}
\ead[label=e3,mark]{qparis@hse.ru}},
\and
\author[C]{\fnms{Nikita} \snm{Zhivotovskiy}
\ead[label=e4]{nikita.zhivotovskii@math.ethz.ch}}

\runauthor{D. Belomestny, L. Iosipoi, Q. Paris, and N. Zhivotovskiy}

\address[A]{Faculty of Mathematics,
Duisburg-Essen University, 
Essen, Germany,
\printead{e1}
}

\address[B]{Faculty of Computer Science,
HSE University,
Moscow, Russia,
\printead{e2,e3}
}

\address[C]{Department of Mathematics, ETH, Z\"{u}rich, Switzerland,
\printead{e4}
}

\end{aug}

\begin{abstract} 
We study the problem of empirical minimization for variance-type functionals over functional classes. Sharp non-asymptotic bounds for the excess variance are derived under mild conditions. In particular, it is  shown that under some restrictions imposed on the functional class  fast convergence rates can be achieved including the optimal non-parametric rates for expressive classes in the non-Donsker regime under some additional assumptions. Our main applications include variance reduction and optimal control.
\end{abstract}
 

\end{frontmatter}

\section{Introduction}\label{sec:introduction}

Empirical Risk Minimization (ERM) plays a central role in statistics and machine learning nowadays.
The goal of learning is usually to find a model which delivers good generalization performance over an underlying distribution of the data. Let $P_{X,Y}$ be a joint distribution of the vector $(X,Y)\in \R^d\times \R$. Given a loss function $\ell:\R\times\R\to\R$, one aims at minimizing the risk $R(h)=\E[\ell(h(X),Y)]$ over a class $\H$ of functions (predictors)  $h:\R^d\to\R$.
Since the distribution $P_{X,Y}$ is usually unknown, one usually replaces  $P_{X,Y}$ by its empirical counterpart and considers 
\begin{eqnarray*}
    h_n\in \argmin_{h\in \H} R_n(h),
\end{eqnarray*}
where  
\begin{eqnarray*}
    R_n(h):=\frac{1}{n}\sum_{i=1}^n \ell(h(X_i),Y_i),
\end{eqnarray*}
and  $\{(X_1,Y_1),\ldots,(X_n,Y_n)\}$ is a sample from $P_{X,Y}$ called ``data'' or ``training set''. ERM covers many popular methods and is widely used in practice. For example, if 
$\H=\{h:\,h(x)=\theta^\top x\ \ \text{for}\ 
\theta\in \R^p\}$ and $\ell(y,p)=(y-p)^2$, then ERM becomes the well-known linear least squares estimator. The celebrated maximum likelihood principle can also be cast as an instance of ERM where the loss function is taken to be the negative log-likelihood function. In turn, ERM itself can be seen as a special case of a more general empirical minimization problem of the form 
\begin{equation}
\label{eq:ustat-opt}
h\in \argmin_{h\in \mathcal{H}} U_n\bigl((h(X_1),Y_1),\ldots, (h(X_n),Y_n)\bigr),
\end{equation}
where $U_n:\R^{2n}\to \R$ is a function which estimates 
$\E[U_n((h(X_1),Y_1),\ldots, (h(X_n),Y_n))]$.
Having the minimal variance among all unbiased estimates
of the mean, see \citep{hoeffding48}, $U$-statistics naturally appear in this context.
Following, for instance, \citep{joly2016robust}, they are defined as
\begin{equation}
\label{eq:ustat}
    U_n=\binom{n}{m}^{-1}\sum_{1\leq i_1<\ldots<i_m\leq n} W\bigl((h(X_{i_1}),Y_{i_1}),\ldots,(h(X_{i_m}),Y_{i_m})\bigr),
\end{equation}
for a fixed positive integer $m \leq n$ and a symmetric measurable function $W:\R^{2m}\to \R$ satisfying $\E|W((h(X_{1}),Y_{1}),\ldots, (h(X_{m}),Y_{m}))|<\infty.$ The summation in \eqref{eq:ustat} is taken over all $m$-element subsets of the set $\{1,\ldots,n\}$.
For the special case $m=2$, and for an appropriate choice of $W$, we recover the empirical variance 
\[
U_n=\frac{1}{n(n-1)}\sum_{1 \leq i<j \leq n}\bigl(Y_i-h(X_i)-Y_j+h(X_j)\bigr)^2,
\]
with $\mathsf{E}[U_n]=\Var[Y-h(X)]$. 
In fact, $U$-statistics generalize common notions of unbiased estimation such as the sample mean, the unbiased sample variance, and the third $k$-statistic which estimates the third cumulant.
After the seminal papers of \citet{halmos46} and \citet{hoeffding48}, it became clear that many meaningful statistics were $U$-statistics. 
We refer to the monograph \citep{korolyuk2013theory}  for a complete exposition of the theory. 
Currently, $U$-statistics play an important role in statistics and probability since they appear in many problems such as 
clustering, ranking, and learning on graphs.
Recently one witnessed a growing interest for concentration properties of $U$-statistics in the context of empirical minimization problems of the form \eqref{eq:ustat-opt}. We refer to \citep*{Clemencon2008},  \citep*{clemencon11},  \citep{clemenccon2014statistical}, \citep{joly2016robust}, \citep{pmlr-v97-clemencon19a}, \citep{pmlr-v99-maurer19a} and references therein. 
\par
In this paper we address the problem of empirical minimization for variance-type functionals, also referred to as Empirical Variance Minimization (EVM). This problem appears in several important applications including variance reduction and optimal control which we discuss in detail in \Cref{sec:applications}. We formulate the problem as follows.
\par
Let $\mathcal X\subseteq \mathbb{R}^d$ and $P$ be a measure on $\mathcal X.$ Furthermore, let  $\H$ be a class  of functions $h:\mathcal X\to \R$ such that $Ph:=\E_P[h(X)]=\mathcal{E}$ for all $h\in \mathcal H$ and some constant $\mathcal{E}.$ The last assumption is a natural convention since the variance as well as the empirical variance are translation invariant. The analysis for a general case follows without changes.
In this paper we study  non-asymptotic properties of the empirical minimizer
\[
	h_n \in \argmin_{h\in\H}\Vn(h)
\]
with 
\begin{align}
	\label{eq:vn}
	\Vn(h)&:=\frac{1}{n-1} \sum_{k=1}^n \bigl(h(X_k) - P_n h\bigr)^2
	= \frac{1}{n(n-1)}\sum\limits _{1\leq i<j\leq n}\bigl(h(X_{i})-h(X_{j})\bigr)^{2},
\end{align}
being the empirical variance based on an i.i.d. sample $X_1,\ldots,X_n$  from $P$ and $P_n h$~---~the empirical mean with respect to the sample $X_1,\ldots,X_n.$
Our goal is to investigate the magnitude of the excess variance  
\[
    \V(h_n)-\V(h^*),
\]
where 
\[
    h^* \in \argmin_{h\in\H} \V(h),
\] 
and $V(h):=\Var_P[h(X)]$ is the true variance of $h(X)$ under $P$. 
For simplicity, we assume that minimizers $h_n$ and $h^*$ exist as all the arguments can easily be adapted by considering approximate minimizers.
We show, for instance, that under suitable conditions on the class  $\H$, the excess variance can be of order up to $O(n^{-1})$ with high probability.  Finally, we consider applications of  these results to problems of variance reduction and optimal control. 
The paper is organized as follows. In \Cref{sec: erm} we present our main results for the EVM problem including a bound on the variance excess under various complexity assumptions on the class \(\H\). \Cref{sec:applications} contains some applications including variance reduction and optimal control. The proofs of the main results are collected in \Cref{sec:proofs}. Some auxiliary results can be found in \Cref{sec:auxres}.
\paragraph*{Notation.}
We use the standard notation for $L^r(P)$-norms, $1 \le r \le \infty$, and denote the set of all functions $h$ with $\| h \|_{L^r(P)}<\infty$ by $L^r(P)$.
We write $\N\bigl(\H, \| \cdot \|_{L^r(P)} ,\eps\bigr)$ for the $\eps$-covering number of 
$\H \subset L^r(P)$, that is, the minimal number of balls of radius $\eps>0$ w.r.t. distance $\| \cdot \|_{L^r(P)} $ needed to cover $\H$.
The natural logarithm of the covering number is called metric entropy of $\H$.
Further, given two functions $g_1, g_2 \in L^r(P)$ such that $g_1 \le g_2$, with probability one, 
the bracket $[g_1,g_2]$ is the set of all functions $h\in L^r(P)$ satisfying $g_1 \le h \le g_2$, with probability one. The size of the bracket $[g_1,g_2]$ is defined as $\| g_1-g_2 \|_{L^r(P)}$. 
The $\eps$-bracketing number $\Nb\bigl(\H, \| \cdot \|_{L^r(P)} ,\eps\bigr)$ of the set 
$\H \subset L^r(P)$ is the minimal number of
brackets of size less or equal to $\eps>0$ necessary to cover $\H$.
The natural logarithm of the bracketing number is called the bracketing  entropy of $\H$. Finally, for any estimator $h_n$ based on the sample $X_1,\ldots, X_n$
the variance $\V(h_n)$ is assumed to be taken conditionally on this sample
and is a $(X_1, \ldots, X_n)$-dependent random variable, that is,
for a new $X \sim P$, independent of $(X_1, \ldots, X_n)$, 
\[
    \V(h_n) := \Var\bigl[h_n(X) \vert X_1,\ldots,X_n\bigr].
\]

\paragraph*{Comparison with previous results and techniques.}
The analysis of empirical minimization of $U$-statistics is explored in several papers. Under the margin assumption, the seminal paper by \citet*{Clemencon2008} provides fast rates (up to $O(n^{-1})$) for $U$-statistics in the context of ranking problems. An important difference with our bounds is that their results take the form of \emph{non-exact} oracle inequalities (see \citep[Theorem 5]{Clemencon2008}), namely, bounds (in our context) on $V(h_n) - 2V(h^*)$. In contrast, our results present exact oracle inequalities, these are the bounds on the excess variance $V(h_n) - V(h^*)$.
Non-exact oracle inequalities are known to be easier to establish, compared to more informative bounds on $V(h_n) - V(h^*)$, and fast rates can be obtained even in the cases where exact oracle inequalities only admit slow rates of convergence \citep{lecue2012}. In particular, non-exact oracle inequalities controlling $V(h_n) - 2V(h^*)$ only require the boundedness of the  loss and exploit neither the convexity of the class nor the margin assumption. In the context of empirical variance minimization, our results provide the first analysis of $U$-statistics that implies the exact oracle inequality with the optimal rate of convergence. More recent results \citep*[Theorem 12]{clemenccon2016scaling} require that the optimal rule belongs to the class in order to provide an exact oracle inequality. In order to simplify the proofs only finite classes are considered in \citep{clemenccon2014statistical} and \citep{joly2016robust}. In contrast, we are considering a larger spectrum of classes: from small parametric classes to expressive classes in the \emph{non-Donsker} regime.
\par
There is also a recent line of research \citep{pmlr-v99-maurer19a} providing bounds for general nonlinear statistics satisfying certain weak interaction assumptions. However, these techniques applied to the special case of $U$-statistics seem to only recover slow rates of convergence.
\par
Finally, the problem of Empirical Variance Minimization was 
announced in \citep{biz2018}. The main difference with this work lies in the fact that the proposed estimators $\tilde{h}_n = \argmin_{h\in\tilde{\H}} \Vn(h)$ were computed by minimization over a finite approximation ($\eps$-net) $\tilde{\H}$ of $\H$, which simplifies the analysis, in the spirit of \citep{devroy95} or \citep{geer00}. 
\par
It has been established that $U$-statistics of order two have concentration properties (for example, the Bernstein inequality) analogous to those of sums of independent random variables \citep{hoeffding63}. 
Our contribution builds upon the observation that a similar phenomenon happens in the localized analysis of empirical minimization of $U$-statistics over expressive classes of functions: the classical rates of parametric and non-parametric regression are also achievable when minimizing certain $U$-statistics. In particular, our analysis avoids referring to the order two Rademacher chaos as in the analysis by \citet{Clemencon2008} and leads to the first sharp oracle inequalities with fast rates of convergence.
\par
On the technical side, this paper presents a self-contained version of the localized analysis of \citet{Kol06} adapted to a more challenging case of $U$-statistics. 
We remark that the original approach of \citeauthor{Kol06} has close similarities to the localized analysis of \citet{massart2007concentration} and of \citet*{bartlett2005}. Our analysis further extends these techniques beyond the standard i.i.d. setting.

\section{Main results}\label{sec: erm}
As we mentioned above, without the loss of generality we consider the class of functions $\H$ such that $Ph =\mathcal{E}$ for all $h\in \mathcal H$ and some constant $\mathcal{E}.$ We work under the following assumptions on the class $\H$. 
\vspace{0.5em}
\begin{enumerate}[label=(A\arabic*),leftmargin=*]
\item\label{A1} There exists $b>0$ such that  $\sup_{h\in\H}\|h\|_{L^{\infty}(P)}\leq b$. \vspace{0.2cm}
\item\label{A2} The class $\H$ is star-shaped around $h^*$, that is, for any $h\in\H$ and any $t\in[0,1]$, 
the function $t h + (1-t)h^*$ belongs to $\H$. 
\end{enumerate}
\vspace{0.5em}
\par
Assumptions~\ref{A1} and \ref{A2} are fairly standard in the analysis based on the empirical process theory (see, for instance, the monographs \citep{Kol11} and  \citep{wainwright_2019}). Note that Assumption~\ref{A2} is clearly satisfied if $\H$ is convex.
\par
Our first result addresses the case where the class $\H$ has either polynomial covering number in $L^2(P_n)$-norm or polynomial bracketing number in $L^2(P)$-norm. These classes are sometimes referred to as  {parametric} or {VC-type} classes (see \citep{barron99}, \citep*{rakhlin2017}). We remark that the results expressed in terms of bracketing numbers is of importance in our application Section \ref{sec:applications}. 
\begin{thm}
\label{th:bigth1}
Suppose that Assumptions \ref{A1} and \ref{A2} hold. 
Suppose also that, for all $0<u\le b$, either
\begin{equation}
\label{eq:polycn}
	\N\bigl(\H,\mn,u\bigr)\le\left(\frac{c}{u}\right)^{\alpha} \ \text{a.s.,}
	\quad\text{or}\quad
	\Nb\bigl(\H,\m,u\bigr)\le\left(\frac{c}{u}\right)^{\alpha}
\end{equation}
for some constants $c>be$ and $\alpha>0$.
Then, for all $n\ge1$ and all $t>0$, with probability at least $1-4e^{-t}$, 
\[
	\V(h_n)-\V(h^*)
	\le A\max\left\{\frac{\log n}{n},\frac{ t}{n}\right\},
\] 
where $A>0$ does not depend on $n, t$ and is explicitly given in the proof up to a multiplicative constant.
\end{thm}
In the above result, the restrictions on the complexity of class $\H$ applies to a fairly large set of situations. However, it  excludes classes $\H$ of non-parametric nature such as H\"{o}lder or Sobolev balls (see, for example, Chapters 2.6 and 2.7 in~\citep{vaart96}). The following result addresses this case and applies to the situation where the
 covering number or bracketing number of $\H$ is exponential in $1/u$.
\begin{thm}
\label{th:bigth2}
Suppose that Assumptions \ref{A1} and \ref{A2} hold. 
Suppose also that, 
for all $0<u\le 2b$, either
\begin{equation}
\label{eq:expcn}
	\log \N\bigl(\H,\mn,u\bigr)\le\left(\frac{c}{u}\right)^{\alpha} \ \text{a.s.,}
	\quad\text{or}\quad
	\log \Nb\bigl(\H,\m,u\bigr)\le\left(\frac{c}{u}\right)^{\alpha}
\end{equation}
for some constants $c>0$ and $\alpha>0$.
Then, for all $n\ge 1$ and all $t>0$, with probability at least $1-4e^{-t}$,
\[
	\V(h_n)-\V(h^*) 
	\le A\max\left\{\vartheta_n,\frac{t}{n}\right\},
	\quad
    \vartheta_n 
    =
    \begin{cases}
        n^{-\frac{2}{2+\alpha}}, &\text{if }\alpha<2,\\
        (\log n) n^{-1/2}, &\text{if }\alpha=2,\\
        n^{-1/\alpha}, &\text{if }\alpha>2,
    \end{cases}
\]
where $A>0$ does not depend on $n, t$ and is explicitly given in the proof up to a multiplicative constant.
\end{thm}
The rates displayed in \Cref{th:bigth1} and \Cref{th:bigth2} are typical for ERM over classes with polynomial and exponential covering numbers respectively. Up to our knowledge, these results are the first fast rates (that is, of order $o(n^{-1/2})$, at least in \Cref{th:bigth1} and in \Cref{th:bigth2} for $\alpha<2$) for empirical variance minimizers. As already mentioned, fast rates in empirical minimization of $U$-statistics have been first obtained by \citet{Clemencon2008} in the context of ranking problems but only in the sense of \emph{non-exact} oracle inequalities.
\par
As was mentioned in the Introduction, the proof of \Cref{th:bigth1} and \Cref{th:bigth2} is based on the localized analysis of \citet{Kol06}. We show that the excess risk $\delta_n:=\V(h_n)-\V(h^*)$ satisfies the so-called fixed point equation
$\delta_n\le \phi_n(\delta_n)$ for some function $\phi_n$ depending on $n$ and on certain ``geometric" and ``complexity" properties of the function class $\H$.
By solving the inequality $\delta \le \phi_n(\delta)$ for $\delta>0$, one can construct a high-probability bound on $\delta_n$ (see Lemma \ref{lem:fixedpoint} which contains a simplification of arguments from \citep{Kol06}; it can be of independent interest). The function $\phi_n$ for general classes $\H$ is given in \Cref{lem:phi}. 
Based on complexity assumptions on $\H$, we upper bound the excess risk $\delta_n$ in \Cref{sec:bigth1} and \Cref{sec:bigth2} for \Cref{th:bigth1} and \Cref{th:bigth2} correspondingly.
Note that the empirical risk $\Vn(h)$ is not a sum of independent random variables $h(X_1),\ldots,h(X_n)$ contrary to classical problems in empirical processes.
Instead, $\Vn(h)$ is a particular case of U-statistics.
This fact complicates the application of the general approach 
from \citep{Kol06}.
\par
The rates $\vartheta_n$ of \Cref{th:bigth2} are well-known in the non-parametric regression literature. In particular, the rate $n^{-2/(2+\alpha)}$ is minimax optimal under various mild assumptions, see \citep*{barron1999risk}. The situation {$\alpha>2$, also known as the \emph{non-Donsker regime} \citep{vaart96}, corresponds to very expressive classes and is more subtle}. It was first shown by \citet{birge1993rates} that empirical minimizers cannot improve the rate $n^{-1/\alpha}$ at least for some special non-parametric classes whereas it is not the case for some specific estimators \citep{barron99}. The possibility of obtaining the rates $n^{-2/(2+\alpha)}$ for all $\alpha > 2$ in the context of non-parametric regression and density estimation attracted a lot of attention (see, for instance, \citep{rakhlin2017}, \citep{pmlr-v65-nikita17a}, \citep{HWCS19}, \citep{kur2019optimality}, and references therein). Following this line of research we provide a similar result for EVM: under the additional assumption that we have access to a well-separated subset of $\H$ with respect to the $L^2(P)$-distance, and that $\H$ contains a function with small variance, we are able to prove that there is a specific variance estimator that achieves the rate $n^{-2/(2+\alpha)}$ for all $\alpha > 0$. This result may be considered as a generalization of the results announced in \citep{biz2018}. Before we proceed, we formulate the following assumption.

\vspace{0.5em}
\begin{enumerate}[label=(A\arabic*), leftmargin=*]\addtocounter{enumi}{2}
\item\label{A3} There is an absolute constant $c_1 > 0$  and a constant $c_2 > 0$, which may depend on the parameters of the problem, such that for any $h \in \H$, 
\[
V(h) - V(h^*) \le c_1\|h - h^*\|_{L^2(P)}^2 + c_2n^{-\frac{2}{2+\alpha}}.
\]
\end{enumerate}
\vspace{0.5em}
This assumption is rather restrictive, though it is always satisfied if the class $\H$ contains a  function with small variance. Indeed, using $(a + b)^2 \le 2a^2 + 2b^2$ and $Ph = Ph^*$, we have 
\begin{equation}
\label{eq:variancedistance}
V(h) - V(h^*) = \|h - Ph\|_{L^2(P)}^2 - \|h^* - Ph\|_{L^2(P)}^2 \le 2\|h - h^*\|_{L^2(P)}^2 + V(h^*).
\end{equation}
Hence, provided that $V(h^*) \le c_2n^{-\frac{2}{2+\alpha}}$ for a fixed constant $c_2>0$, we see that \ref{A3} holds with $c_1 = 2.$ Assuming that the variance $V(h^*)$ is decreasing as the sample size $n$ grows is quite natural in some situations. Indeed, due to the non-asymptotic nature of our results, one may choose $\H$ according to the value of the sample size. When the sample size $n$ is getting larger, it is reasonable to consider a more expressive class $\H$. As a result, the minimum variance in the class can decrease. We are now ready to state our last excess variance bound which can be seen as a strengthening of \Cref{th:bigth2} under additional assumptions.
\begin{thm}
\label{th:bigth3}
Suppose that Assumptions \ref{A1}, \ref{A2}, \ref{A3} hold.
Suppose also that, for all $0<u\le b$,
\[
	\log \N\bigl(\H,\m,u\bigr)\le\left(\frac{c}{u}\right)^{\alpha}
\]
for some constants $c > 0$ and $\alpha>0$.
Then there is an estimator $\tilde{h}_n$ such that, for all $n\ge 1$ and all $t>0$, with probability at least $1-4e^{-t}$,
\[
	\V(\tilde{h}_n)-\V(h^*) 
	\le A\max\left\{n^{-\frac{2}{2+\alpha}}, \frac{t}{n}\right\},
\] 
where $A>0$ does not depend on $n, t$, but depends on $c_2$ in Assumption \ref{A3} and is explicitly given in the proof. 
\end{thm}
\begin{remark*}
Observe that in \Cref{th:bigth3} the assumption $\log \Nb\bigl(\H,\m,u\bigr)\le(c/u)^{\alpha}$ is replaced by the weaker assumption $\log \N\bigl(\H,\m,u\bigr)\le(c/u)^{\alpha}$.
\end{remark*}
\begin{remark*}
It follows from the proof of \Cref{th:bigth3} that in the special case where $V(h^*)=0$, that is, $\H$ contains a constant function, Assumption \ref{A2} is not required and Assumption \ref{A3} is automatically satisfied due to \eqref{eq:variancedistance}.
\end{remark*}
Let us now present the estimator of \Cref{th:bigth3}. 
We fix
\begin{equation}
\label{def:eta}
   \eps = b^{\frac{2}{2+\alpha}}\,n^{-\frac{1}{2+\alpha}},
\end{equation}
and consider the set $\H_{\eps} \subseteq \H$ to be a minimal $\eps$-net of $\H$ with respect to the $L^2(P)$ distance. Define
\begin{equation}
    \label{eq:sieveestimator}
\tilde{h}_n = \argmin\limits_{h \in \H_{\eps}}V_n(h).
\end{equation}
The estimators of the form \eqref{eq:sieveestimator} are analyzed in statistical literature and
are referred to as the skeleton or sieve estimators. The known results imply the bounds for
procedures used in density estimation, non-parametric regression and classification (see \citep*{wong95}, \citep*{devroy95}, and \citep*{geer00}). Observe that the construction of the set $\H_{\eps}$ depends on the $L^2(P)$ structure of $\H$ and therefore, requires some prior knowledge on $P$ or an access to additional unlabelled data points which allows to estimate the $L^2(P)$ distances between functions in $\H$. 

\section{Applications}
\label{sec:applications}
In this section we provide several applications of our general results.
\subsection{Variance reduction}
Suppose that we wish to compute $Pf$, 
where $f:\X\to\R$ is a function from $L^2(P)$. 
A natural estimate for $Pf$ is the Monte Carlo estimate
$P_nf$.
This estimate is unbiased, that is, $\E [P_nf] = Pf$, and has variance $\Var[P_nf]=\V(f)/n$. From the point of view of applications, an important problem is to improve accuracy of the Monte Carlo estimate $P_nf$,
which in turn means reducing its variance.
This can be achieved by increasing the sample size $n$ but in some cases this solution may not be practical.
Another approach to reduce the variance is to decrease the value $\V(f)$ 
by considering a new Monte Carlo estimate with the same mean but lower variance. 
Such techniques are called variance reduction methods, see
\citep{christian1999monte}, \citep{rubinstein2016simulation}, and \citep{glasserman2013monte} for  an introduction to this field.
The method of control variates is one of the few generic variance reduction methods.
The idea behind this method is to choose a class $\G \subset L^2(P)$ of functions such that any
$g\in\G$ satisfies $Pg=0$ (such functions are called control variates), and to
consider a new Monte Carlo estimate 
\[
	P_n\bigl(f-g^*\bigr) = \frac1n \sum_{k=1}^n \bigl(f(X_k)-g^*(X_k)\bigr),
\] 
where 
\begin{equation*}
	g^* \in \argmin\limits_{g\in\G}\V\bigl(f-g\bigr).
\end{equation*} 
If the class $\G$ of control variates is properly chosen, this approach may significantly reduce the variance, 
that is, $\V(f-g^*) \ll \V(f)$.
Since the true variance $\V(f-g)$ for all $g\in\G$ is not assumed to be known
and its computation is a more difficult problem than the initial one, a natural approach to pick a
control variate is to minimize the empirical variance  
\begin{equation}
	g_n \in \argmin\limits_{g\in\G}\Vn\bigl(f-g\bigr).
	\label{eq:vrgn}
\end{equation} 
Consequently, since the true variance is replaced by its empirical counterpart, the variance of $f-g_n$ can be decomposed as
\begin{equation*}
	\V(f - {g}_n) =  \underbrace{\V(f - {g}_n)  - \inf_{g\in\G} \V(f - g)}_{\text{stochastic error}}  
	+ \underbrace{\inf_{g\in\G} \V(f - g)}_{\text{approximation error}}.
\end{equation*} 
So there is a natural trade-off between the stochastic and approximation errors. Namely, as the set $\G$ increases the approximation error decreases while the estimation
error typically increases. In the particular case where the function $g^* = f-Pf$ belongs to $\G$ the approximation error is equal to zero, hence the variance of $f-g_n$ includes only the stochastic error. 
\par
In order to bound the stochastic error, \Cref{th:bigth1}, \Cref{th:bigth2}, or \Cref{th:bigth3} may be used with 
$\H = \{ f-g:\ g\in\G\}$. 
These results show that, under some technical assumptions, the stochastic error is of order $O(n^{-1})$ when $\H$ has a polynomial covering number as in \eqref{eq:polycn} and
ranges from $O(n^{-2/2+\alpha})$ for $\alpha<2$ to $O(n^{-1/\alpha})$ for $\alpha>2$ when $\H$ has an exponential covering number as in \eqref{eq:expcn} (note that covering numbers of $\H$ and $\G$ are the same).
Similar results were previously announced in \citep*{biz2018} but only for the sieve estimator $\tilde{g}_n$ of the form \eqref{eq:sieveestimator} under more restrictive assumptions. 
Their results also imply that the stochastic error is of order up to $O(n^{-1})$. 
In \citep*{bimns2020, bimns2020-2} the authors study a more sophisticated
setting where $\tilde{g}_n$ is estimated based on dependent samples generated by MCMC algorithms. In this case the stochastic error for the sieve estimator is bounded by $O(n^{-1/2})$. 
\par
The effect of control variates on the Monte Carlo estimate $P_nf$ can be illustrated by constructing confidence intervals before and after the variance reduction. 
Namely, using a  sample of length $N$ and the central limit theorem, the following asymptotic confidence interval for $Pf$ holds
\begin{equation}
\label{eq:ci-stand}
	P_Nf \pm q\,\sqrt\frac{\V(f)}{N},
\end{equation}
where $q$ is a quantile of the standard normal distribution. 
If we now apply the EVM variance reduction approach with a class \(\G\) having polynomial covering number and construct \(g_n\) using a new independent sample of the length \(n\), then the confidence interval becomes 
\begin{equation*}
	P_N(f-g_n) \pm q\,\sqrt\frac{\V(f-g_n)}{N}=P_Nf \pm q\,\sqrt\frac{\inf_{g\in\G} \V(f - g)}{N}+O\left(\sqrt{\frac{1}{nN}}\right),
\end{equation*}
which can be significantly tighter than \eqref{eq:ci-stand} if \(n\) is large.
Loosely speaking, this approach may reduce the length of the asymptotic confidence interval from order $O(N^{-1/2})$ to $O(n^{-1/2}N^{-1/2})$, provided that the class $\G$ is chosen so that
$\inf_{g\in\G} \V(f - g) $ is small enough
(this actually means that the class of control variates $\G$ possesses nice approximation properties).
Of course, it comes at a price --- this procedure requires computational resources to choose a good control variate (that is, to solve the optimization problem \eqref{eq:vrgn}).
\par
Furthermore, let us examine separately the case of Stein control variates (see \citep{assaraf1999}, \citep{mira2013zero}, and \citep{oates2017control}), where the class $\G$ is constructed by substitution
of various smooth functions $\phi: \X \to \R^d$ into
\begin{equation*}
	g_\phi= \langle\phi,\nabla \log\pi\rangle + \divergence(\phi),
\end{equation*}
where $\pi$ is the density of the distribution $P$, $\langle\cdot,\cdot\rangle$ denotes the standard inner product in $\R^d$, and $\divergence(\phi)$ stands for the divergence of $\phi$. 
Under rather mild conditions on $\pi$ and $\phi$, integration by parts implies 
that $Pg_\phi=0$ (see \citep[Propositions~1 and 2]{mira2013zero}). 
In what follows, $W^{s,p}(\lambda)$ denotes the Sobolev space with respect to the Lebesgue measure $\lambda$, that is,
$
	W^{s,p}(\lambda) = \bigl \{ u \in L^p(\lambda) : D^{\alpha}u \in L^p(\lambda) \ \text{for all} \ |\alpha| \leqslant s \bigr \},
$
where $\alpha = (\alpha_1,\ldots,\alpha_d)$ is a multi-index with $|\alpha|=\alpha_1+\ldots+\alpha_d$ and 
$D^{\alpha}$ is the differential operator of the form
$D^{\alpha} = \partial^{|\alpha|}/\partial x_1^{\alpha_1}\ldots\partial x_d^{\alpha_d}$ 
(here partial derivatives are understood in the generalized sense of distributions).
We write $W^{s,p}_R(\lambda)$ for the class of functions $h\in W^{s,p}(\lambda)$
supported on a closed ball $B_R=\{x\in \X: |x|\leq R\}$.
The Sobolev norm for $u\in W^{s,p}(\lambda)$ or $u\in W^{s,p}_R(\lambda)$ is defined as
$\|u\|_{W^{s,p}(\lambda)} = \sum\nolimits_{|\alpha|\le s}\|D^{\alpha}u\|_{L^p(\lambda)}$. By a slight abuse of notation, we continue to write $\phi\in W_R^{s,p}(\lambda)$ for  functions $\phi:\, \X \to \R^d$ meaning that any coordinate of $\phi$ belongs to $W_R^{s,p}(\lambda)$. 
Given a class $\Phi$
of smooth functions $\phi:\X\to\R^d$, define the corresponding minimizers
\[  
    \phi_n \in \argmin\limits_{\phi\in\Phi}\Vn\bigl(f-g_\phi\bigr)
    \quad\text{and}\quad
    \phi^* \in \argmin\limits_{\phi\in\Phi}\V\bigl(f-g_\phi\bigr).
\]
\begin{prop}\label{th:corvr}
Fix some \(R>0\) and let $\Phi\subset W^{s,2}_R(\lambda)$ be a class of $s$-smooth functions $\phi:\X\to\R^d$ with $s-d/2>0$.
Let also $\pi\in W^{1,2}(\lambda)$ be a probability density function satisfying
$1/L\leq \pi(x)\leq L$ for some $L>0$ and all $x\in B_R$. 
If $\H = \{ f-g_\phi:\, \phi\in\Phi\}$ satisfies \ref{A1} and \ref{A2}, then for all $n\ge 1$ and $t>0$,
with probability at least $1-4e^{-t}$,
\[
	\V(f-g_{\phi_n})-\V(f-g_{\phi^*}) 
	\le A\max\left\{\vartheta_n,\frac{b^2t}{n}\right\},
	\quad
	\vartheta_n 
    =
    \begin{cases}
        n^{-\frac{2}{2+d^2/s}}, &d^2/s<2,\\
        (\log n) n^{-1/2}, &d^2/s=2,\\
        n^{-s/d^2}, &d^2/s>2,
    \end{cases}
\]
where $A>0$ is a constant not depending on $n$.
\end{prop}

\begin{remark*}
It is worth mentioning that in \Cref{th:corvr} we consider the case where $\Phi$ has large complexity. If, for example, covering number of $\Phi$ is polynomial, 
the rates may be improved by applying \Cref{th:bigth1} instead of \Cref{th:bigth2} (assumption $\Phi\subset W^{s,2}_R(\lambda)$ corresponds to exponential covering number). Furthermore, if one has an access to a minimal $\eps$-net of $\Phi$, \Cref{th:bigth3} can be applied.
In this case $\tilde{h}_n$ will not depend on any properties of $f$ other than its values observed on the sample.
\end{remark*}

\subsection{Optimal control}
We consider the optimal control of a discrete-time Markov process with a finite time horizon \(T\). On
a filtered measurable probability space $(\Omega,\mathcal{F})$ with $\mathcal{F}=\left(
\mathcal{F}_{r}\right)_{r=0,1,\ldots,T}$, $T\in\mathbb{Z}_{+}$, we define an
adapted control process $\mathbf{a}:\,\Omega\times\{0,\ldots,T-1\}\rightarrow A,$
control for short, where $(A,\mathcal{B})$ is a measurable space. We assume a given set of admissible controls which is denoted by $\mathcal{A}.$ Given a control $\mathbf{a}=(a_{0},a_{1},...,a_{T-1})\in\mathcal{A}$, we consider a
\textit{controlled Markov} process $X$ valued in some measurable space
$(S,\mathcal{S})$ and defined on a probability space $(\Omega,\mathcal{F},\P ^{\mathbf{a}})$ with $X_{0}=x_{0}$ a.s. and transition kernel of the
type
\[
P^{a_{r}}(x,dy)=\P ^{\mathbf{a}}(X_{r+1}\in dy\ |\ X_{r}=x),\quad 0\leq r<T.
\]
So, it is assumed that the distribution of $X_{r+1}$ conditional
on $\mathcal{F}_{r}$ is governed by a (one-step) transition kernel $P^{a_{r}}(X_{r},dy)$ which is in turn controlled by $a_{r}.$ In this way to each admissible
control function, one associates a random variable and a measure. One can therefore also associate to the admissible control \(\mathbf{a}\) the mean \(\E^{\mathbf{a}}\) of the alternative random variable with respect to the alternative measure \(\P ^{\mathbf{a}}\). 
This procedure defines a continuous function from the set of admissible control functions to the real numbers. We refer to this function as the reward of the underlying optimal control problem. The desired unknown value that we wish to estimate is equal to the global maximum 
of the reward with respect to all admissible controls. In this way we may consider
the general optimal control problem of the form:
\begin{equation}
Y_{0}^{\ast}=\sup_{\mathbf{a}\in\mathcal{A}}\E^{\mathbf{a}}\left[ \, \sum
_{r=0}^{T-1}f_{r}(X_{r},a_{r})+g_T(X_T)\right] 
\label{eq:opt}
\end{equation}
for given functions  $f_{r}$, $r=0,\ldots, T-1,$ and $g_T.$
Introduce the process 
\begin{equation}
Y_{r}^{\ast}=\sup_{\mathbf{a}\in\mathcal{A}_{r}}\E^{\mathbf{a}}\left[ \,\sum_{s=r}^{T-1}f_{s}
(X_{s},a_{s})+g_{T}(X_{T}) \,\bigg\vert\, \mathcal{F}_{r}\right]
,\quad
0\leq r\leq T ,
\label{genopt1}
\end{equation}
with  $\mathcal{A}_{r}$ being the set of
all admissible controls $\mathbf{a}:\Omega\times\{r,\ldots,T-1\}\to A.$ Then
there exists a vector $h^{\ast}=(h_{0}^{\ast},\ldots
,h_{T}^{\ast})$ of measurable functions on $S,$ such that $Y_{j}^{\ast}
=h_{j}^{\ast}(X_{j})$ and $h^{\ast}$ satisfies the dynamic programming principle:
\begin{equation*}
h_{r}^{\ast}(x)    =\left(  \mathcal{L}h^{\ast}\right)
_{r}(x)  ,\quad0\leq r<T,\quad 
h_{T}^{\ast}(x) =g_{T}(x), 
\end{equation*}
where $\mathcal{L}$ is a Bellman-type operator
defined by
\[
(\mathcal{L}h)  _{r}(x)=\sup_{a\in A\,}\left[  f_{r}
(x,a)+P^{a}h_{r+1}(x)\right]
\]
and
\[
P^{a}h_{r+1}(x):=\int P^{a}(x,dy)h_{r+1}(y).
\]
We now assume that there exists a reference measure $\P ^{\ast}$ equivalent to
$\P ^{\mathbf{a}},$  such that
\[
P^{a}(x,dy)=\varphi(x,y,a)P^{\ast}(x,dy),\quad a\in A,
\]
with $P^{\ast}(x,dy):=\P ^{\ast}(X_{r+1}\in dy\ |\ X_{r}=x)$ and the function
$\varphi(x,y,a)$ satisfying $\varphi\geq 0$ and 
\[
\int P^{\ast}(x,dy)\varphi
(x,y,a)\equiv1.
\]
\par
Denote by \(C_b(S)\) the set of continuous bounded functions on \(S.\)
As shown in \citep{rogers2007pathwise} the value of the problem \eqref{eq:opt} 
in a dual  form can be expressed as an infimum over a family of  martingales  \(M_j=\left(  \mathcal{L}h\right)
_{j}(X_{j})-h_{j}(X_{j}),\)  \(j=1,\ldots,T,\) \(h\in C_b(S)\)  of an expectation, under measure \(P^{\ast}\), of a pathwise supremum  adjusted by a weighted  sum of \(M_j\). 
\begin{thm*}[\citeauthor{rogers2007pathwise}]
\label{Dualrep} 
Let $Y^{*}_{r} $ be a solution of the optimal control
problem \eqref{genopt1}, then the following representation holds
\begin{equation}
\label{eq:dual-opt}
Y_{0}^{\ast}  =\inf_{h\in C_b^{T+1}(S)}\mathsf{E}^{\ast}\left[  h_{0}(X_{0})+ \sum_{j=0}^{T-1}W_{j}\left( \left(  \mathcal{L}h\right)
_{j}(X_{j})-h_{j}(X_{j})\right) \right]  ,
\end{equation}
where 
\[
W_{j}=\sup_{\mathbf{a}\in\mathcal{A}}\left[ \, {\displaystyle\prod\limits_{l=0}^{j-1}}
\varphi(X_{l},X_{l+1},a_{l}) \right] 
\]
and \(\mathsf{E}^{\ast}\) stands for expectation under \(P^{\ast}.\)
\end{thm*}
   As shown in \citep{belomestny2018dual} the optimization problem \eqref{eq:dual-opt} can be alternatively formulated as a problem of variance minimization 
\begin{equation}
\label{eq:dual-var}
Y_{0}^{\ast}  =\inf_{h\in C_b^{T+1}(S)}\Var\left [h_{0}(X_{0})+\sum_{j=0}^{T-1}W_{j}\left( \left(  \mathcal{L}h\right)
_{j}(X_{j})-h_{j}(X_{j})\right) \right],
\end{equation}
provided the set \(A\) of control values is finite. To solve this problem, one can use Monte Carlo approach and cast the problem \eqref{eq:dual-var} into the empirical variance minimization 
\[
h_n\in \arginf_{h\in\mathcal{H}}\left\{\frac{1}{n(n-1)}\sum_{1\leq k<l\leq n} \bigl(\phi_h(X^{(k)}_0,\ldots,X^{(k)}_{T})-\phi_h(X^{(l)}_0,\ldots,X^{(l)}_{T})\bigr)^2\right\}, 
\]
where   \(\mathcal{D}_n=\bigl\{(X^{(k)}_{0},\ldots,X^{(k)}_{T}), \,k=1,\ldots,n\bigr\}\) 
is a set of \(n\) independent trajectories of the chain \((X_j)_{j\geq 0},\) $\H\subset C_b^{T+1}(S)$, and  
\[
\phi_h(X_0,\ldots,X_{T})=h_{0}(X_{0})+\sum_{j=0}^{T-1}W_{j}\left( \left(  \mathcal{L}h\right)
_{j}(X_{j})-h_{j}(X_{j})\right). 
\]
\Cref{th:bigth1} implies that, for all $n\ge1$ and all $t>0$, with probability at least $1-4e^{-t}$, 
\[
	\Var\bigl[\phi_{h_n}(X_0,\ldots,X_{T})|\mathcal{D}_n\bigr]-\inf_{h\in\H} \Var\bigl[\phi_{h}(X_0,\ldots,X_{T})\bigr]
	\le A\max\left\{\frac{\log n}{n},\frac{ t}{n}\right\}
\] 
for some $A>0$ not depending on $n, t$,
provided that the class of functions $\Phi=\{\phi_h: \,h\in \H\}$  satisfies the assumptions \ref{A1}, \ref{A2} and 
\[
\log \Nb\bigl(\Phi,\|\cdot\|_{L_2(P^*\times \ldots \times P^*)},u\bigr)\le\left(\frac{c}{u}\right)^{2}.
\] 
Note that the latter condition is, for example, fulfilled if $\Phi$ is a ball in Sobolev space 
$W^{s,2}(S^{T+1})$ for $s$ large enough. This condition in turn can be translated into a smoothness assumption imposed on the class $\H$ and the transition densities $\varphi(x,y,a),$ $a\in A,$ see the proof of Proposition~\ref{th:corvr}. 

\section{Proofs}
\label{sec:proofs}

\subsection{Outline of the proofs of \texorpdfstring{\Cref{th:bigth1} and \Cref{th:bigth2}}{Lg}}

Below, we present an outline of the proof of \Cref{th:bigth1} and \Cref{th:bigth2} based on 
technical results proved in \Cref{sec:auxres}. 
The proof consists in an adaptation of the general strategy for bounding empirical risk minimizers 
devised in~\citep{Kol06}. The main difference being that we deal with the quadratic functionals $\V(h)$ and $\Vn(h)$ 
instead of $Ph$ and $P_nh$ respectively. Essentially, we prove that the strategy extends nicely to our setting using Hoeffding decomposition for $U$-Statistics.
Let 
\[
	\delta_n:=\V(h_n)-\V(h^*).
\] 
Notice that by the definition of $h_n$, we get
\begin{align*}
	\delta_n &\le \bigl(\V(h_n)-\V(h^*)\bigr)-\bigl(\Vn(h_n)-\Vn(h^*)\bigr)\\
	&= P\bigl((h_n-\e)^2-(h^*-\e)^2\bigr)-\bigl(\Vn(h_n)-\Vn(h^*)\bigr),
\end{align*}
where we recall that $\e = Ph$ for all $h\in\H$.
Observe that, for all $h\in\H$, $\Vn(h)$ may be rewritten as
\begin{align*}
	\Vn(h) &= \frac{1}{2n(n-1)} \sum_{i \ne j} \bigl( h(X_i) - h(X_j) \bigr)^2\\
	&= P_n(h-\e)^2-\frac{1}{n(n-1)} \sum_{i \ne j} (h(X_i)-\e)(h(X_j)-\e).
\end{align*} 
This expression is a simple instance of the well-known 
Hoeffding decomposition for $U$-statistics (due to \citet{hoeffding48} and \citet{hoeffding61}). 
Combining this representation with the above inequality yields
\[
	\delta_n \le T_n(h_n)+W_n(h_n),
\]
where 
\[
	T_n(h):=(P-P_n)\bigl((h-\e)^2-(h^*-\e)^2\bigr),
\]
and 
\[
	W_n(h)=w_n(h) - w_n(h^*)
	\quad\text{with}\quad
	w_n(h)=\frac{1}{n(n-1)} \sum_{i \ne j} (h(X_i)-\e)(h(X_j)-\e).
\]
As a result, we arrive at the central observation that the excess risk $\delta_n$ satisfies
\begin{align*}
    \delta_n\le \phi_n(\delta_n),
\end{align*}
where, for any $\delta>0$,
\begin{equation}
\label{eq:hdeltadef}
	\phi_n(\delta):=\sup_{h\in\H(\delta)}\bigl\{T_n(h)+W_n(h)\bigr\}
	\quad\text{and}\quad 
	\mathcal H(\delta)=\bigl\{h\in\H:\ \V(h) - \V(h^*) \le \delta\bigr\}.
\end{equation}

The main technical result we will invoke, in order to bound $\delta_n$, is the following lemma. 
This result follows from a combination of arguments presented in Theorem~4.1, Corollary~4.1, and Theorem~4.3 in \citep{Kol11}. 
We simplify these arguments and, for the sake of completeness, provide the proof in \Cref{sec:lemfixedpoint}.
\begin{lem}
\label{lem:fixedpoint}
Let $\{\phi(\delta):\,\delta\ge 0\}$ be non-negative random variables (indexed by all deterministic $\delta\ge 0$) 
such that, almost surely, $\phi(\delta)\le \phi(\delta')$ if $\delta\le\delta'$.
Let $\{\beta(\delta,t):\, \delta\ge 0, t\ge0\}$, be (deterministic) real numbers such that 
\begin{equation}
\label{lem:fixedpoint:e1}
	\P\bigl(\phi(\delta)\ge \beta(\delta,t)\bigr)\le e^{-t}.
\end{equation}
Finally, let $\hat \delta$ be a nonnegative random variable, a priori upper bounded by a constant $\bar\delta>0$, and such that, almost surely, 
\[
	\hat \delta\le \phi(\hat \delta).
\]
Then defining, for all $t\ge0$,
\begin{equation}
\label{lem:fixedpoint:e2}
	\beta(t):=\inf\left\{\tau>0:\sup_{\delta\ge\tau}\frac{\beta\left(\delta,\tfrac{t\delta}{\tau}\right)}{\delta}\le \frac{1}{2}\right\},
\end{equation}
we obtain, for all $t\ge0$, 
\[
	\P\bigl(\hat \delta\ge \beta(t)\bigr)\le 2e^{-t}.
\]
\end{lem}
\par 
The result of \Cref{lem:fixedpoint} is a version of the fixed point argument which is widely applied in the statistical literature to prove the rates of convergence faster than $O(n^{-1/2})$. 
For a deeper discussion on this argument
we refer to the papers \citep{birge1993rates}, \citep{barron1999risk}, \citep{bartlett2005} and the monographs \citep{geer00}, \citep{vaart96}, \citep{massart2007concentration}, \citep{wainwright_2019}.
\par
According to \Cref{lem:fixedpoint}, it remains to bound $\phi_n(\delta)$ with high probability for any fixed $\delta>0$. The next lemma reduces the problem to that of bounding the expected suprema of empirical processes. Below, for every $h\in\mathcal H$, we denote
\begin{equation}
\label{eq:ellh}
   \ell(h):=(h-\e)^2-(h^*-\e)^2. 
\end{equation}
\begin{lem}
\label{lem:phi}
Suppose that Assumptions \ref{A1} and \ref{A2} hold. 
Then, for every $\delta>0$ and any $t>0$,
\[
	\P\bigl(\phi_n(\delta)\ge \beta_n(\delta,t)\bigr)\le 2e^{-t},
\]
where, 
\[
	\beta_n(\delta,t)
	:=
	A\left( \E\sup_{h\in\H(\delta)}(P - P_n)\ell(h) + \left(\E\sup_{h\in\H(\delta)}|(P-P_n)h|\right)^2 +b\sqrt{\frac{t\delta}{n}}+\frac{b^2 (1+t)}{n}\right),
\]
$\ell(h)$ is defined by \eqref{eq:ellh}, and $A>0$ is a universal constant.
\end{lem}
\begin{proof}
Let us write
\[
	\phi_n(\delta)\le \phi^{(1)}_n(\delta)+\phi^{(2)}_n(\delta),
\]
where
\[
	\phi^{(1)}_n(\delta):=\sup_{h\in\H(\delta)}T_n(h)
	\quad\text{and}\quad 
	\phi^{(2)}_n(\delta):=\sup_{h\in\H(\delta)}W_n(h).
\]
\noindent\textbf{\emph{Bound on $\phi^{(1)}_n(\delta)$}}.
Under \ref{A1}, we have $|\ell(h)| \le 4b^2$ for all $h\in\H$. 
It then follows from a version of Talagarand's inequality due to Bousquet (see \Cref{lem:bousquet} in \Cref{sec:bousquet}) that, with probability at least $1-e^{-t}$,
\[
	\phi^{(1)}_n(\delta)\le \E \phi^{(1)}_n(\delta)+\sqrt{\frac{2t}{n}\bigl(\sigma^2(\delta)+16b^2\,\E \phi^{(1)}_n(\delta)\bigr)}+\frac{8b^2 t}{3n},
\]
where \[\sigma^2(\delta) = \sup_{h\in\mathcal H(\delta)} P\ell(h)^2.\] 
Using basic inequalities $\sqrt{u+v}\le \sqrt{u}+\sqrt{v}$ and $2\sqrt{uv}\le u+v$ for positive numbers $u$ and $v$, we further deduce that, with probability at least $1-e^{-t}$,
\begin{equation}
	\label{eq:phi1:1}
	\phi^{(1)}_n(\delta)\le 2\E \phi^{(1)}_n(\delta)+\sigma(\delta)\sqrt{\frac{2t}{n}}+\frac{32b^2 t}{3n}.
\end{equation}
Now let us upper bound the term $\sigma^2(\delta)$. Observe that under \ref{A1}, for every $h\in\H$,
\[
	|\ell(h)|=|(h-h^*)(h+h^*-2\e)|\le 4b|h-h^*|,
\]
so that 
\[
	P\bigl(\ell(h)\bigr)^2\le 16b^2P(h-h^*)^2.
\]
Then, using the simple identity $(h-h^*)^2=2\ell(h)-4\ell((h+h^*)/2)$, we deduce that
\begin{equation}
\label{eq:bernsteinassumption}
    P(h-h^*)^2=2P\ell(h) - 4 P\ell\left(\frac{h+h^*}{2}\right)
    \le 2P\ell(h),
\end{equation}
where the inequality follows by definition of $h^*$ and the fact that,
since $\mathcal H$ is star-shaped around $h^*$ according to \ref{A2}, the function $(h+h^*)/2$ belongs to $\mathcal H$. Hence we have proven that, for every $h\in\mathcal H$,
\begin{equation*}
	P\bigl(\ell(h)\bigr)^2 \le 32b^2P\ell(h).
\end{equation*}
This inequality
is an instance of the Bernstein assumption \citep{bartlett2006empirical}.
By the definition \eqref{eq:hdeltadef} of $\H(\delta)$, it follows that
\begin{equation}
	\label{eq:phi1:2}
	\sigma^2(\delta) \le 32b^2\delta. 
\end{equation}
Combining \eqref{eq:phi1:1} and \eqref{eq:phi1:2}, we deduce that 
\begin{equation}
    \label{eq:phi1:3}
    \phi^{(1)}_n(\delta)\le 2\E\sup_{h\in\mathcal H(\delta)}(P-P_n)\ell(h)+8b\sqrt{\frac{t\delta}{n}}+\frac{32b^2 t}{3n}.
\end{equation}
\noindent\textbf{\emph{Bound on $\phi^{(2)}_n(\delta)$}}.
Notice that, for all $h\in\mathcal H$,
\begin{align*}
w_n(h)&=
    \frac{1}{n(n-1)} \sum_{i=1}^n  \sum_{j=1}^n  (h(X_i)-\e)(h(X_j)-\e)
    -
    \frac{1}{n(n-1)} \sum_{i=1}^n (h(X_i)-\e)^2\\
    &=
    \frac{n}{n-1}\bigl(P_n(h-\e)\bigr)^2-\frac{1}{n-1}P_n(h-\e)^2\\
    &=
    \frac{n}{n-1}\bigl((P-P_n)h\bigr)^2-\frac{1}{n-1}P_n(h-\e)^2.
\end{align*}
Omitting negative terms and using \ref{A1} we get, for all $h\in\mathcal H$ and all $n\ge 2$, 
\begin{align*}
   	W_n(h)
   	&\le\frac{n}{n-1}\bigl((P-P_n)h\bigr)^2+\frac{1}{n-1}P_n(h^*-\e)^2\\
   	&\le 2\bigl((P-P_n)h\bigr)^2+\frac{8b^2}{n}.
\end{align*}
As a result, for any $\delta>0$, we have
\begin{align*}
    \phi^{(2)}_n(\delta)
    	&\le 2\left(\sup_{h\in\H(\delta)}| (P-P_n)h|\right)^2+\frac{8b^2}{n}\\
    	&\le 4\left(\sup_{h\in\H(\delta)}| (P-P_n)h|-\E\sup_{h\in\H(\delta)}| (P-P_n)h|\right)^2
    	+4\left(\E\sup_{h\in\H(\delta)}| (P-P_n)h|\right)^2+\frac{8b^2}{n}.
\end{align*}
A classical application of the bounded differences inequality \citep[Theorem 6.2]{boucheron2013concentration} implies that, with probability at least $1-e^{-t}$,
\begin{equation*}
	\left(\sup_{h\in\H(\delta)}| (P-P_n)h|-\E\sup_{h\in\H(\delta)}| (P-P_n)h|\right)^2 \leq \frac{2b^2t}{n}.
\end{equation*}
Hence, we deduce from the above that, for every $\delta>0$ and every $t>0$,
\begin{equation}
	\label{eq:phi2}
	\phi^{(2)}_n(\delta)\le 4\left(\E\sup_{h\in\H(\delta)}|(P-P_n)h|\right)^2+ \frac{8b^2(1+t)}{n},
\end{equation}
with probability at least $1-e^{-t}$. 
Combining \eqref{eq:phi1:3} and \eqref{eq:phi2} we obtain the desired conclusion.
\end{proof}

The rest of the proofs of \Cref{th:bigth1} and \Cref{th:bigth2} is given in \Cref{sec:bigth1} and \Cref{sec:bigth2} correspondingly.
There, depending on the assumption on the entropy of the class $\H$, we explicitly compute $\beta_n(\delta,t)$ defined in \Cref{lem:phi} and substitute it in \Cref{lem:fixedpoint} to obtain the final bound on the excess risk $\delta_n$.

\subsection{End of the proof of \texorpdfstring{\Cref{th:bigth1}}{Lg}} 
\label{sec:bigth1}
Throughout the proof, let $A>0$ be a universal constant whose value may change from line to line. 
It follows from \Cref{lem:phi} that, for every $\delta\ge0$ and every $t>0$,
\[
	\P\bigl(\phi_n(\delta)\ge \beta_n(\delta,t)\bigr)\le 2e^{-t},
\]
where
\[
	\beta_n(\delta,t):=A\left( \E\sup_{h\in\H(\delta)}(P - P_n)\ell(h) + \left(\E\sup_{h\in\H(\delta)}|(P-P_n)h|\right)^2 +b\sqrt{\frac{t\delta}{n}}+\frac{b^2 (1+t)}{n}\right).
\]
Let us now bound the expected suprema of empirical processes given in $\beta_n(\delta,t)$ under the assumptions of \Cref{th:bigth1} on either the covering or 
bracketing numbers.  
Under \ref{A1}, we have $|\ell(h)|\leq4b^2$ for any $h\in\H$. 
Moreover, for any $h_1,h_2\in\H$, 
\[
	|\ell(h_1) - \ell(h_2)|\leq|(h_1-h_2)(h_1+h_2-2\e)|\leq4b|h_1-h_2|,
\]
and hence 
\[
	\N\bigl(\L(\delta),\mn,u\bigr)\leq\N\biggl(\H(\delta),\mn,\frac{u}{4b}\biggr).
\]
The same result holds also for the bracketing numbers. 
From \eqref{eq:phi1:2} we also see that 
\[
    \sup_{h\in\H(\delta)} P\bigl(\ell(h)\bigr)^2 \leq \sqrt{32b^2\delta}.
\]
Combining these facts and taking  $\sigma = \sqrt{32b^2\delta}$ 
in either \Cref{lem:boundrad1} for covering numbers or \Cref{lem:boundrad3} 
for the bracketing numbers (both results are presented in \Cref{sec:auxres}), we get
\begin{align*}
	\E\sup_{h\in\mathcal H(\delta)}(P - P_n)\ell(h)
	&\leq 
	A \left(\sqrt{\frac{32\alpha b^2 \delta}{n}\log\left(\frac{4c}{\sqrt{32\delta}}\right)} + \frac{\alpha b}{n}\log\left(\frac{4c}{\sqrt{32\delta}}\right)\right)\\
	&\leq
	A \left(\sqrt{\frac{\alpha b^2 \delta}{n}\log\left(\frac{c^2A}{\delta}\right)} + \frac{\alpha b}{n}\log\left(\frac{c^2A}{\delta}\right)\right).
\end{align*}
Now let us turn to the second summand in $\beta_n(\delta,t)$. 
Letting $\bar{\H}(\delta):=\H(\delta)\cup(-\H(\delta))$, we obtain
$\N\bigl(\bar{\H}(\delta),\mn,u\bigr)\leq2\N\bigl(\H(\delta),\mn,u\bigr)$
and
\[
	 \E\sup_{h\in\H(\delta)}|(P-P_n)h| = 
	 \E\sup_{h\in\bar{\H}(\delta)}(P-P_n)h.
\]
Using either \Cref{lem:boundrad1} for the covering number assumption 
or \Cref{lem:boundrad3} for the bracketing number assumption, we get with $\sigma = b$ 
that
\[
	\E\sup_{h\in\H(\delta)}|(P-P_n)h| \leq 
	A \left(\sqrt{\frac{\alpha b^2}{n}\log\left(\frac{c}{b}\right)} + \frac{\alpha b}{n}\log\left(\frac{c}{b}\right)\right).
\]
Hence, whenever $n \ge \alpha \log(c/b)$,
\[
	\left(\E\sup_{h\in\H(\delta)}|(P-P_n)h|\right)^2  
	\leq A \frac{\alpha b^2}{n}\log\left(\frac{c}{b}\right).
\]
Combining these we conclude
\[
	\beta_n(\delta,t) 
	\leq A\left( 
	\sqrt{\frac{\alpha b^2 \delta}{n}\log\left(\frac{c^2A}{\delta}\right)} 
	+ \frac{\alpha b^2}{n}\log\left(\frac{c^2A}{\delta}\right) 
	+ b\sqrt{\frac{t\delta}{n}}+\frac{b^2 (1+t)}{n}
	\right).
\]
\Cref{lem:fixedpoint} now implies that 
\[
	\P\bigl(\delta_n\ge \beta_n(t)\bigr)\le 4e^{-t}
	\quad\text{with}\quad
	\beta_n(t):=\inf\left\{\tau>0:\sup_{\delta\ge \tau}\frac{\beta_n\left(\delta,\tfrac{t\delta}{\tau}\right)}{\delta}\le \frac{1}{2}\right\}.
\]
It remains only to compute an upper bound for $\beta_n(t)$. To that aim, it may be easily checked that, for any $\tau>0$, 
\[
	\sup_{\delta\ge \tau}\frac{\beta_n\left(\delta,\tfrac{t\delta}{\tau}\right)}{\delta}
	=
	A\left(
	\sqrt{\frac{\alpha b^2}{\tau n}\log\left(\frac{c^2A}{\tau}\right)}
	+\frac{\alpha b^2}{\tau n}\log\left(\frac{c^2A}{\tau}\right)
	+b\sqrt{\frac{t}{\tau n}}+\frac{b^2(1+t)}{\tau n}
	\right).
\]
Next, notice that for any non-increasing functions $y_1,\dots,y_k$, we have
\[
	\inf\{\tau>0:y_{1}(\tau)+\dots+y_{k}(\tau)\le 1\}\le \max_{1\le i\le k}\inf\left\{\tau>0:y_{i}(\tau)\le \frac{1}{k}\right\}.
\]
As a result, we deduce that 
\[\beta_n(t)\le \max_{1\le i\le 4}\beta^{(i)}_n(t),\]
where
\begin{align*}
    \beta^{(1)}_n(t)&= \inf\left\{\tau>0: A\sqrt{\frac{\alpha b^2}{\tau n}\log\left(\frac{c^2A}{\tau}\right)}\le \frac{1}{8}\right\},\\
    \beta^{(2)}_n(t)&=\inf\left\{\tau>0:A\frac{\alpha b^2}{\tau n}\log\left(\frac{c^2A}{\tau}\right)\le \frac{1}{8}\right\},\\
    \beta^{(3)}_n(t)&=\inf\left\{\tau>0:Ab\sqrt{\frac{t}{\tau n}}\le \frac{1}{8}\right\},\\
    \beta^{(4)}_n(t)&=\inf\left\{\tau>0:A\frac{b^2(1+t)}{\tau n}\le \frac{1}{8}\right\}.
\end{align*}
To upper bound $\beta^{(1)}_n(t)$ and $\beta^{(2)}_n(t)$, we use the fact that for all $u>0$ and $v>0$,
\[
	\inf\left\{\tau>0:\frac{u}{n\tau}\log\frac{v}{\tau}\le 1\right\}\le \frac{u \log n}{n},
\]
as soon as $\log n\ge v/u$ (this fact can be checked by a substitution of the given bound). 
For $\beta^{(3)}_n(t)$ and $\beta^{(4)}_n(t)$, we compute the infimums
directly. We finally obtain
\begin{align*}
     \beta^{(1)}_n(t) \leq\frac{64A^2\alpha b^2 \log n}{n},
     \quad
     \beta^{(2)}_n(t)\leq \frac{8A \alpha b^2 \log n}{n},
     \quad
     \beta^{(3)}_n(t)= \frac{64 A^2 b^2 t}{n},
     \quad
     \beta^{(4)}_n(t)=\frac{8 A b^2 (1+t)}{n},
\end{align*}
which finishes the proof for large enough $n$,
that is, when $\log n \ge Ac^2 / (\alpha b^2)$.
For smaller $n$ the desired formula obviously holds for large $A>0$. 

\subsection{End of the proof of \texorpdfstring{\Cref{th:bigth2}}{Lg}}
\label{sec:bigth2}
As before, throughout the proof, let $A>0$ be a universal constant 
whose value may change from line to line. 
It follows from \Cref{lem:phi} that, for every $\delta\ge0$ and every $t>0$,
\[
	\P\bigl(\phi_n(\delta)\ge \beta_n(\delta,t)\bigr)\le 2e^{-t},
\]
where, for a universal constant $A>0$,
\[
	\beta_n(\delta,t):=A\left( \E\sup_{h\in\H(\delta)}(P - P_n)\ell(h) + \left(\E\sup_{h\in\H(\delta)}|(P-P_n)h|\right)^2 +b\sqrt{\frac{t\delta}{n}}+\frac{b^2 (1+t)}{n}\right).
\]
Let us now bound the expected suprema of empirical processes given in $\beta_n(\delta,t)$ under the assumptions of \Cref{th:bigth2} on either the metric or 
bracketing entropy. 
Consider the following cases of values $\alpha$.\\ \\
\noindent\textbf{\emph{Case $\alpha<2$}}.
Here the proof repeats in many ways the proof of \Cref{th:bigth1}, 
so some details are omitted. 
Taking $\sigma = \sqrt{32b^2\delta}$ 
in either \Cref{lem:boundrad2} for the metric entropy assumption or 
 \Cref{lem:boundrad4} for the bracketing entropy assumption, we get 
\begin{align*}
	\E\sup_{h\in\H(\delta)}(P - P_n)\ell(h)
	&\leq 
	A \left(\sqrt{\frac{b^2\delta}{(2-\alpha)^2\,n}\left(\frac{c}{\sqrt{\delta}}\right)^{\alpha}} + \frac{b}{(2-\alpha)^2\,n}\left(\frac{c}{\sqrt{\delta}}\right)^{\alpha}\right).
\end{align*}
Again depending on the entropy assumption we use, 
either \Cref{lem:boundrad2} or \Cref{lem:boundrad4}
with $\sigma = b$ yields 
\[
	\E\sup_{h\in\H(\delta)}|(P-P_n)h|
	\leq 
	A \left(\sqrt{\frac{b^2}{(2-\alpha)^2\,n}\left(\frac{c}{b}\right)^{\alpha}} + \frac{b}{(2-\alpha)^2\,n}\left(\frac{c}{b}\right)^{\alpha}\right).
\]
Hence, whenever $n \ge (2-\alpha)^{-2} (c/b)^{\alpha}$, 
\[
	\left(\E\sup_{h\in\H(\delta)}|(P-P_n)h|\right)^2  
	\leq
	A \, {\frac{b^2}{(2-\alpha)^2\,n}\left(\frac{c}{b}\right)^{\alpha}}.
\]
Combining these bounds we conclude
\[
	\beta_n(\delta,t) 
	\leq A\left( 
	\sqrt{\frac{b^2\delta}{(2-\alpha)^2\,n}\left(\frac{c}{\sqrt{\delta}}\right)^{\alpha}} + \frac{b^2}{(2-\alpha)^2\,n}\left(\frac{c}{\sqrt{\delta}}\right)^{\alpha} 
	+ b\sqrt{\frac{t\delta}{n}}+\frac{b^2 (1+t)}{n}
	\right).
\]
\Cref{lem:fixedpoint} now implies that 
\[
	\P\bigl(\delta_n\ge \beta_n(t)\bigr)\le 4e^{-t}
	\quad\text{with}\quad
	\beta_n(t):=\inf\left\{\tau>0:\sup_{\delta\ge \tau}\frac{\beta_n\left(\delta,\tfrac{t\delta}{\tau}\right)}{\delta}\le \frac{1}{2}\right\}.
\]
It remains only to compute an upper bound for $\beta_n(t)$. To that aim, it may be easily checked that, for any $\tau>0$, 
\[
	\sup_{\delta\ge \tau}\frac{\beta_n\left(\delta,\tfrac{t\delta}{\tau}\right)}{\delta}
	=
	A\left( 
	\sqrt{\frac{b^2}{(2-\alpha)^2\,\tau n}\left(\frac{c}{\sqrt{\tau}}\right)^{\alpha}} 
	+ \frac{b^2}{(2-\alpha)^2\,\tau n}\left(\frac{c}{\sqrt{\tau}}\right)^{\alpha}
	+b\sqrt{\frac{t}{\tau n}}
	+\frac{b^2 (1+t)}{\tau n}
	\right).
\]
As a result, we deduce that 
\[\beta_n(t)\le \max_{1\le i\le 4}\beta^{(i)}_n(t),\]
where 
\begin{align*}
    \beta^{(1)}_n(t)&= \inf\left\{\tau>0: A \sqrt{\frac{b^2}{(2-\alpha)^2\,\tau n}\left(\frac{c}{\sqrt{\tau}}\right)^{\alpha}} \le \frac{1}{8}\right\} = \left(\frac{64A^2b^2c^\alpha}{(2-\alpha)^2\,n}\right)^{\frac{2}{2+\alpha}},\\
    \beta^{(2)}_n(t)&=\inf\left\{\tau>0:A\frac{b^2}{(2-\alpha)^2\,\tau n}\left(\frac{c}{\sqrt{\tau}}\right)^{\alpha}\le \frac{1}{8}\right\}= \left(\frac{8Ab^2c^\alpha}{(2-\alpha)^2\,n}\right)^{\frac{2}{2+\alpha}},\\
    \beta^{(3)}_n(t)&=\inf\left\{\tau>0:Ab\sqrt{\frac{t}{\tau n}}\le \frac{1}{8}\right\}=\frac{64 A^2 b^2 t}{n},\\
    \beta^{(4)}_n(t)&=\inf\left\{\tau>0:A\frac{b^2 (1+t)}{\tau n}\le \frac{1}{8}\right\}=\frac{8 A b^2 (1+t)}{n}.
\end{align*}
This completes the proof for $\alpha < 2$.\\ \\
\noindent\textbf{\emph{Case $\alpha=2$}}.
Repeated application of either \Cref{lem:boundrad2} for the metric entropy assumption or \Cref{lem:boundrad4} for the bracketing entropy assumption allows us to write,
whenever $n\ge (b \log n /c)^2$, that
\[
	\beta_n(\delta,t) 
	\leq A\left( 
	\frac{cb\log n}{\sqrt{n}} + \frac{b}{n}
	+ b\sqrt{\frac{t\delta}{n}}+\frac{b^2 (1+t)}{n}
	\right).
\]
\Cref{lem:fixedpoint} now implies that 
\[
	\P\bigl(\delta_n\ge \beta_n(t)\bigr)\le 4e^{-t},
\]
where
\begin{align*}
	\beta_n(t)
	\leq 
	A\max\left\{ \frac{cb\log n}{\sqrt{n}}, \frac{b}{n}, \frac{b^2t}{n}, \frac{b^2(1+t)}{n}\right\}.   
\end{align*}
This finishes the proof for $\alpha = 2$.\\ \\
\noindent\textbf{\emph{Case $\alpha>2$}}.
Similarly, repeated application of either \Cref{lem:boundrad2} for the metric entropy assumption or \Cref{lem:boundrad4} for the bracketing entropy assumption enables us to write,
whenever 
\[
n\ge \max\{b^\alpha, c^{\alpha/2}b^{-\alpha(\alpha-1)}(\alpha-2)^{-\alpha}\},
\]
that
\[
	\beta_n(\delta,t) 
	\leq A\left( 
	\frac{b}{n^{1/\alpha}} + \frac{c^{\alpha/2}b}{(\alpha-2)\,n^{1/\alpha}}
	+ b\sqrt{\frac{t\delta}{n}}+\frac{b^2 (1+t)}{n}
	\right).
\]
\Cref{lem:fixedpoint} now implies that 
\[
	\P\bigl(\delta_n\ge \beta_n(t)\bigr)\le 4e^{-t},
\]
where
\begin{align*}
	\beta_n(t)
	\leq 
	A \max\left\{ \frac{b}{n^{1/\alpha}}, \frac{c^{\alpha/2}b}{(\alpha-2)\,n^{1/\alpha}}, \frac{b^2t}{n}, \frac{b^2(1+t)}{n}\right\}.   
\end{align*}
This completes the proof of \Cref{th:bigth2}.

\subsection{Proof of \texorpdfstring{\Cref{th:bigth3}}{Lg}}
The main technical result we will invoke is the Bernstein-type inequality for $U$-statistics (see inequality A.1 in \citep{Clemencon2008}). By definition,
we can write for any $h\in\H$ that
\[
    \Vn(h) - \Vn(h^*)  = \frac{1}{n(n-1)} \sum_{1\leq i, j \leq n}\frac{\bigl( h(X_i) - h(X_j) \bigr)^2 - \bigl( h^*(X_i) - h^*(X_j) \bigr)^2}{2} .
\]
Using \ref{A1}, the kernel of this $U$-statistic satisfies for any $1\leq i, j \leq n$,
\begin{equation*}
\frac{\bigl|\bigl(h(X_i) - h(X_j)\bigr)^2 - \bigl(h^*(X_i) - h^*(X_j)\bigr)^2\bigr|}{2}\le 2b^2.
\end{equation*}
Now the Bernstein inequality immediately implies that, for all $\delta \ge 0$,
\begin{equation}
\label{eq:bernsteinineq}
    \P\Bigl(\bigl| \Vn(h) - \Vn(h^*) - \bigl(\V(h) - \V(h^*)\bigr)\bigr| \ge  \delta \Bigr) \le 2\exp\left(-\frac{n\delta^2}{\sigma^2 + 4b^2\delta /3}\right), 
\end{equation}
where $\sigma^2 = \Var\bigl[(h(X_i) - h(X_j))^2 - (h^*(X_i) - h^*(X_j))^2\bigr]$.
A simple computation gives
\begin{align*}
    \sigma^2 
    &\le \E\Bigl[\bigl((h(X_i) - h(X_j))^2 - (h^*(X_i) - h^*(X_j))^2\bigr)^2\Bigr]
    \\
    &\le 16 b^2 \, \E\Bigl[\bigl(h(X_i) - h(X_j) - h^*(X_i) + h^*(X_j)\bigr)^2\Bigr]
    \\
    &\le 32 b^2 \, P(h - h^*)^2.
\end{align*}
Using the fact that $P(h - h^*)^2 \leq 2 (\V(h) - \V(h^*))$, see
\eqref{eq:bernsteinassumption} in the proof of \Cref{lem:phi},
we get
\[
    \sigma^2 \leq 32 b^2 \, P(h - h^*)^2 \leq 64b^2 \bigl(\V(h) - \V(h^*)\bigr).
\]
The main idea of the proof is to consider concentration of $\Vn(h) - \Vn(h^*)$ not around the mean $\V(h) - \V(h^*)$, but around
$c(\V(h) - \V(h^*))$ for some $c>1$. This idea is often used to obtain fast rates, see, for example, \citep{bartlett2005} and \citep{Kol06}.
Combining the last bound for $\sigma^2$ with \eqref{eq:bernsteinineq} we obtain for all $\delta>0$,
\begin{align}
    &\P\Bigl(\Vn(h) - \Vn(h^*) - 2\bigl(\V(h) - \V(h^*)\bigr)\ge \delta\Bigr) 
    \nonumber\\&\qquad\qquad\le
    \P\Bigl(\Vn(h) + \Vn(h^*) - \bigl(\V(h) - \V(h^*)\bigr)  \ge \delta + \V(h) - \V(h^*)\Bigr)
    \nonumber\\&\qquad\qquad\le 
    2\exp\left(-\frac{n(\delta + V(h) - V(h^*))^2}{64b^2(V(h) - V(h^*)) + 4b^2(\delta + V(h) - V(h^*))/3}\right)
    \nonumber\\
    &\qquad\qquad\le 2\exp\left(-\frac{n(\delta + V(h) - V(h^*))}{(64 + 4/3)b^2}\right)
    \nonumber\\
    &\qquad\qquad\le 2\exp\left(-\frac{n\delta}{66b^2}\right),
    \label{eq:bernsteinineq1}
\end{align}
where the last inequality follows from the fact that $V(h) - V(h^*) \ge 0$.
Similarly,
\begin{align}
    &\P\Bigl(\V(h) - \V(h^*) - 2\bigl(\Vn(h) - \Vn(h^*)\bigr)\ge \delta\Bigr) 
    \nonumber\\&\qquad\qquad\le
    \P\left(\V(h) + \V(h^*) - \bigl(\Vn(h) - \Vn(h^*)\bigr)  \ge \frac{\delta + \V(h) - \V(h^*)}{2}\right)
    \nonumber\\&\qquad\qquad\le 
    2\exp\left(-\frac{n\delta}{260b^2}\right).
    \label{eq:bernsteinineq2}
\end{align}
Having disposed of this preliminary step, we can now return to the proof.
Fix any $\eps>0$ and let $\H_{\eps}$ be a  minimal $\eps$-net  of $\H$ 
with respect to $L_2(P)$. Choose any $\tilde{h}_{n}\in\argmin_{h\in\H_\eps} \Vn(h)$ and any $h^*_{\eps} \in \H_{\eps}$  such that $\|h^*_{\eps} - h^*\|_{L_2(P)} \le \eps$. 
We have the following decomposition
\begin{align}
    \V(\tilde{h}_{n}) - \V(h^*) &\le \V(\tilde{h}_{n}) - \V(h^*) - 2\bigl(\Vn(\tilde{h}_{n}) - \Vn(h^*_{\eps})\bigr)
    \nonumber\\
    &= \V(\tilde{h}_{n}) - \V(h^*) - 2\bigl(\Vn(\tilde{h}_{n}) - \Vn(h^*)\bigr) + 2\bigl(\Vn(h^*_{\eps}) - \Vn(h^*)\bigr)
    \nonumber\\
    &\le \sup\limits_{h \in \H_{\eps}}\bigl\{\V(h) - \V(h^*) - 2(\Vn(h) - \Vn(h^*))\bigr\} + 2\bigl(\Vn(h^*_{\eps}) - \Vn(h^*)\bigr).
    \label{eq:vdecomp}
\end{align}
Let us bound all the terms in \eqref{eq:vdecomp} separately. 
First, we get by \eqref{eq:bernsteinineq2} and the union bound,
\[
    \P\biggl(\sup_{h\in\H_{\eps}}\bigl\{\V(h) - \V(h^*) - 2\bigl(\Vn(h) - \Vn(h^*)\bigr)\bigr\}\ge \delta\biggr) 
    \leq 
    2 \, \N(\H,\m,\eps) \exp\left(-\frac{n\delta}{260b^2}\right).
\]
Choosing $\delta = 260 b^2n^{-1}(\N(\H,\m,\eps) + t) $ and using the assumption on the covering number of $\H$, we obtain that, with probability at least $1-2e^{-t}$,
\[
    \sup_{h\in\H_{\eps}}\bigl\{\V(h) - \V(h^*) - 2\bigl(\Vn(h) - \Vn(h^*)\bigr)\bigr\}
    \leq 
    \frac{260 b^2}{n}\biggl(\frac{c}{\eps}\biggr)^{\alpha} + \frac{260 b^2t}{n}.
\]
Further, for the second term in \eqref{eq:vdecomp} we have
\[
    2\bigl(\Vn(h^*_{\eps}) - \Vn(h^*)\bigr) = 
    2\bigl(\Vn(h^*_{\eps}) - \Vn(h^*) \bigr) 
    -4\bigl(\V(h^*_{\eps}) - \V(h^*) \bigr) 
    +4\bigl(\V(h^*_{\eps}) - \V(h^*) \bigr).
\]
It follows from \eqref{eq:bernsteinineq1} that, with probability at least $1-2e^{-t}$,
\[
    2\bigl(\Vn(h^*_{\eps}) - \Vn(h^*) \bigr) 
    -4\bigl(\V(h^*_{\eps}) - \V(h^*) \bigr) 
    \leq 
    \frac{132b^2t}{n}.
\]
Assumption \ref{A3} implies that
\[
V(h^*_{\eps}) - V(h^*) \le c_1 \eps^2 + c_2n^{-\frac{2}{2+\alpha}}.
\]
Combining all the bounds, we conclude that, with probability at least $1-4e^{-t}$,
\[
    \V(\tilde{h}_{n}) - \V(h^*) \leq \frac{260b^2}{n}\biggl(\frac{c}{\eps}\biggr)^{\alpha} + \frac{392b^2t}{n} + 4c_1 \eps^2 + 4c_2n^{-\frac{2}{2+\alpha}}.
\]
The proof is finished by taking $\eps = b^{2/(2+\alpha)}n^{-1/(2+\alpha)}$.

\subsection{Proof of \texorpdfstring{\Cref{th:corvr}}{Lg}}
With notation $\H = \{ f-g_{\phi},\, \phi\in\Phi\}$ 
and $\G = \{ g_{\phi},\, \phi\in\Phi \}$, 
it obviously holds that
\[
    \Nb(\H,\m,u) = \Nb(\G,\m,u).
\]
For any two functions $\phi_1,\phi_2\in\Phi$,
we have
\begin{align*}
    \| g_{\phi_1}-g_{\phi_2} \|_{L^2(P)}
    &\leq
    \bigl\| \langle \phi_1 - \phi_2, \nabla\log\pi  \rangle \bigr\|_{L^2(P)} + \bigl\|\divergence(\phi_1 - \phi_2)\bigr\|_{L^2(P)}.
\end{align*}
Let us denote the coordinates of $\phi_k:\X \to \R^d$ by
$\phi^i_k$, $i=1,\ldots,d$, $k=1,2$.
By our assumptions and the Cauchy-Schwarz inequality, the first term can be bounded as
\begin{align*}
    \bigl\| \langle \phi_1 - \phi_2, \nabla\log\pi \rangle \bigr\|_{L^2(P)}
    &\leq
    \sqrt{L} \biggl( \int_{B_R} | \phi_1(x) - \phi_2(x) |^2 |\nabla\pi(x)|^2\, \rmd x \biggr)^{1/2}\\
    &\leq
    M\sqrt{L} \sum_{i=1}^d 
    \biggl( \int_{B_R} | \phi^i_1(x) - \phi^i_2(x) |^2  \rmd x \biggr)^{1/2},
\end{align*}
provided that $\|\pi\|_{W^{1,2}(\lambda)}\leq M$ for some $M>0$. Further,
using the Minkowski inequality, 
\begin{align*}
    \bigl\|\divergence(\phi_1 - \phi_2)\bigr\|_{L^2(P)}
    &=
    \Biggl( \int_{B_R} \biggl(\sum_{i=1}^d \frac{\partial}{\partial x_i}( \phi^i_1(x) - \phi^i_2(x))\biggr)^2 \pi(x)\, \rmd x \Biggr)^{1/2}\\
    &\leq
    \sqrt{L}\, \sum_{i=1}^d
     \biggl( \int_{B_R} \Bigl( \frac{\partial}{\partial x_i}( \phi^i_1(x) - \phi^i_2(x))\Bigr)^2 \rmd x \biggr)^{1/2}.
\end{align*}
Combining these bounds, we conclude
\begin{align*}
    \| g_{\phi_1}-g_{\phi_2} \|_{L^2(P)}
    &\leq
    \sqrt{L} (M + 1)
    \sum_{i=1}^d
    \bigl\| \phi^i_1 - \phi^i_2\|_{W^{1,2}(\lambda)}.
\end{align*}
By the arithmetic properties of the bracketing entropy, for all $u>0$,
\[
\Nb(\G,\m,u) \leq \prod_{i=1}^d \Nb\left(\Phi^i,\|\cdot\|_{W^{1,2}(\lambda)},\frac{u}{\sqrt{L} (M + 1) d}\right),
\]
where $\Phi^i$ denotes the class generated by $i$-th coordinate of $\phi\in\Phi$, that is, $\Phi^i = \{ \phi^{i}, \phi\in\Phi\}$. Using the classic results for the bracketing entropy of Sobolev spaces defined on $B_R$ (see, for example, \citep[Corollary~4 and Section~3.3.4]{nickl2007bracketing}), we get
\[
\Nb(\G,\m,u) \leq \prod_{i=1}^d
c \left( \frac{ \sqrt{L} (M + 1) d}{u} \right)^{d/s}
\leq
\left( \frac{c^{s/d}\sqrt{L}(M + 1)d}{u} \right)^{d^2/s}
\]
for some $c>0$ not depending on $u$.
The result now follows from \Cref{th:bigth2}.

\paragraph*{Acknowledgments.} 
This article was prepared within the framework of the HSE University
Basic Research Program. Results
of Section 3 were obtained by Denis Belomestny and Leonid Iosipoi under support of the RSF grant 19-71-30020 (HSE University).
Nikita Zhivotovskiy is funded in part by ETH Foundations of Data Science (ETH-FDS).

\bibliographystyle{plainnat}
\bibliography{penal_vr}
\newpage
\begin{appendix}
\section{Supplementary Material}
\label{sec:auxres}
In this section we have compiled some standard facts on empirical processes. 
Namely, we start with Bousquet's concentration inequality given in \Cref{sec:bousquet}. 
Then we provide bounds on the expectation of suprema of empirical process based on either covering number or bracketing
number assumptions in \Cref{sec:metricbounds} and \Cref{sec:bracketbounds} respectively. 
Throughout the section, $\H$ stands for a class of functions $h:\X\to \R$.

\subsection{Bousquet's concentration inequality}
\label{sec:bousquet}
We start with the well-known concentration result, known as Bousquet's form of Talagrand's inequality for empirical processes. 
It involves a notion of variance of the empirical process 
\[
	\sigma_{\H}:=\sup_{h\in\H} \sqrt{Ph^2},
\]
which plays a crucial role in many modern proof techniques involving the local behaviour of the supremum 
of empirical process. The proof of the following lemma can be found in \citep{Bou02} or \citep{GinNic16}. 

\begin{lem}
\label{lem:bousquet}
Suppose that all functions in $\H$ are $[a, b]$-valued, for some $a<b$.
Then, for all $n\ge 1$ and all $t>0$, 
\[
	\sup_{h\in\H} (P-P_{n})h
	\le \E \sup_{h\in\H} (P-P_{n})h + \sqrt{\frac{2t}{n}\Bigl(\sigma^{2}_{\H}+2(b-a)\,\E \sup_{h\in\H} (P-P_{n})h\Bigr)}+\frac{(b-a)t}{3n},
\]
with probability larger than $1-e^{-t}$.
\end{lem}

\subsection{Bounds on expected supremum of empirical process using metric entropy}\label{sec:metricbounds}
In this section we derive explicit bounds for the expected supremum of empirical process
under specific assumptions on the metric entropy of $\H$. We let $\rv_1,\dots,\rv_n$ be a sequence of i.i.d. random signs, that is, 
$\P(\rv_i=-1)=\P(\rv_i=1)=1/2$, independent from the sample $X_1,\ldots,X_n$. We denote
\[
	R_n(h):= \frac{1}{n} \sum_{i=1}^n\rv_i h(X_i),
\]
the Rademacher process indexed by $\H$.
The next result is a generalized version of Dudley's entropy bound; its proof can be found, for instance, 
in \citep[Lemma~A.3]{srebro2010}.

\begin{lem}
\label{lem:chain}
Let $X_1,\ldots,X_n$ be a sample associated to the empirical distribution $P_n$ and suppose that $\H\subset L^2(P_n)$.
Then, for all $\eps>0$,
\[
	\E\left[\sup_{h\in\H}R_n(h) \Big\vert X_1,\ldots,X_n \right]
	\le 
	4\eps+\frac{12}{\sqrt n}\int_{\eps}^{\sigma_n}\sqrt{\log \N\bigl(\H,\mn,u\bigr)}\rmd u,
\]
where $\sigma_n:=\sup_{h\in\H}\sqrt{P_nh^2}$.
\end{lem}
\par
The next two lemmas apply \Cref{lem:chain} to derive the upper bounds on the expected supremum of the empirical process under the polynomial and exponential covering number assumptions. 
They consist in a simple version of Theorem 3.12 in \citep{Kol11} and can be also found in, for example, \citep{GinKol06}. 
We start with the case of a class $\H$ with polynomial covering number.

\begin{lem}
\label{lem:boundrad1}
Let $\H$ be a class of $[-b,b]$-valued functions for some $b>0$. Suppose that, for all $0< u\le b$,  
\begin{equation*}
	\N\bigl(\H,\mn,u\bigr) \le \left(\frac{c}{u}\right)^{\alpha}
	\ \text{almost surely},
\end{equation*}
for some constants $c>b\rme$ and $\alpha>0$. Then, for any  $\sigma \in [\sigma_{\H},b]$ and all $n\ge 1$,
\[
	\E\sup_{h\in\H} (P-P_{n})h
	\le
	A \left(\sqrt{\frac{\alpha\sigma^2}{n}\log\left(\frac{c}{\sigma}\right)} + \frac{\alpha b}{n}\log\left(\frac{c}{\sigma}\right)\right),
\]
where $A>0$ is an absolute constant computed explicitly in the proof. 
\end{lem}
\begin{proof} 
To simplify the notation throughout the proof, we denote
\[
	R:= \E\sup_{h\in\H}R_n(h).
\] 
From the symmetrization principle, we get 
\begin{align*}
    \E\sup_{h\in\H} (P-P_{n})h
    \le 2R.
\end{align*}
Now it remains to bound $R$.
Taking $\eps=0$ in \Cref{lem:chain}, yields
\begin{equation*}
R \le \frac{12}{\sqrt n}\,\E \left[\int_{0}^{\sigma_n}\sqrt{\alpha\log \left(\frac{c}{u}\right)}\rmd u\right],
\end{equation*}
where we have used the assumption on the covering numbers of $\H$ and the fact that, almost surely, $\sigma_n\leq b$.
Since the function 
$$\eps \mapsto \int_{0}^{\eps}\sqrt{\alpha\log \left(\frac{c}{u}\right)} \rmd u,$$
is concave on $(0,c]$, we deduce from Jensen's inequality that
\begin{equation*}
	R \le \frac{12}{\sqrt n}\,\int_{0}^{\E[\sigma_n]}\sqrt{\alpha\log \left(\frac{c}{u}\right)} \rmd u.
\end{equation*}
Now let us bound the value of $\E[\sigma_n]$.
Using once again Jensen's inequality, we obtain
\begin{align*} 
	\E[\sigma_n]&=\E\left[\sup\nolimits_{h\in\H}\sqrt{P_nh^2}\right]
	\le \sqrt{\E\sup\nolimits_{h\in\H}P_nh^2}.
\end{align*}
The symmetrization and contraction principles imply
\begin{align*} 
	\E[\sigma_n]
	&\le \sqrt{\sigma^2+\E\sup\nolimits_{h\in\H}(P_n-P)h^2}\\
	&\le\sqrt{\sigma^{2}+2\E\sup\nolimits_{h\in\H} R_n(h^2)}\\
	&\le\sqrt{\sigma^{2}+4bR}\\
	&\le\sigma+2\sqrt{bR}.
\end{align*}
Since $\sigma_n\le b$ almost surely, we deduce that 
\[
	\E[\sigma_n]\le B:=\min\Bigl\{b,\sigma+2\sqrt{bR}\Bigr\}.
\]
Consequently, after a change of variable we get
\begin{align} 
	R
	&\le 12\sqrt{\frac{\alpha}{n}}\int_{0}^{B} \sqrt{\log \left(\frac{c}{u}\right)} \rmd u
	=12c\sqrt{\frac{\alpha}{n}}\int_{c/B}^{+\infty} u^{-2}\sqrt{\log(u)} \rmd u.
\label{entrop:e6}
\end{align}
Observe that, by assumptions, $c/B\ge e$. Using integration by parts, it follows that
\[
	\sup_{t\ge e}\frac{\int_{t}^{+\infty} u^{-2}\sqrt{\log(u)}\rmd u}{t^{-1}\sqrt{\log(t)}}
	\leq 
	1 + \sup_{t\ge e}\frac{\frac12\int_{t}^{+\infty} u^{-2}(\log(u))^{-3/2} \rmd u}{t^{-1}\sqrt{\log(t)}}
	\leq
	\frac32.
\]
Therefore, we have
\begin{align}
	\int_{c/B}^{+\infty} u^{-2}\sqrt{\log(cu)}\,{\rm d}u &\le \frac32 \frac{B}{c}\sqrt{\log  \left(\frac{c}{B}\right)}
	\le \frac32 \frac{B}{c}\sqrt{\log \left(\frac{c}{\sigma}\right)},
\label{entrop:e7}
\end{align}
where, in the last inequality, we have used that $B\ge\sigma$. Now combining \eqref{entrop:e6} and \eqref{entrop:e7}, we deduce that 
\begin{equation}
	R \le 18(\sigma+2\sqrt{bR})\sqrt{\frac{\alpha}{n}\log\left(\frac{c}{\sigma}\right)}.
\nonumber
\end{equation}
Since $R\le x+y\sqrt R$ implies $R\le 2x + 2y^{2}$ for any $x,y>0$, we obtain 
\begin{equation*}
	R \le 36\sqrt{\frac{\alpha\sigma^2}{n}\log \left( \frac{c}{\sigma}\right)} + 2592\,\frac{\alpha b}{n}\log \left( \frac{c}{\sigma}\right),
\end{equation*}
which completes the proof.
\end{proof}

Next we study the case of classes $\H$ for which the covering number grows exponentially. 
The important observation here is that the square root of the entropy might not be integrable around zero anymore.
\begin{lem}
\label{lem:boundrad2}
Let $\H$ be a class of $[-b,b]$-valued functions for some $b>0$. Suppose that, for all $0< u\le b$,  
\begin{equation*}
	\log\N\bigl(\H,\mn,u\bigr) \le \left(\frac{c}{u}\right)^{\alpha},
\end{equation*}
for some constants $c>0$ and $\alpha>0$. 
Then for some absolute constant $A>0$ and all $n\ge1$
the following holds.\\
\begin{enumerate}
	\item If $\alpha < 2$, then, for any  $\sigma \in [\sigma_{\H},b]$, 
	\[
		\E\sup_{h\in\H} (P-P_{n})h
		\le
		A\left(
		\frac{\sigma}{(2-\alpha)\sqrt{n}}\left(\frac{c}{\sigma}\right)^{\alpha/2}
		+ \frac{b}{(2-\alpha)^2\,n}\left(\frac{c}{\sigma}\right)^{\alpha}
		\right).
	\]
	\item  If $\alpha = 2$, then
	\[
		\E\sup_{h\in\H} (P-P_{n})h
		\le
		A\left(
		\frac{c\log n}{\sqrt n} + \frac{b}{n} 
		\right).
	\]	
	\item  If $\alpha > 2$, then 
	\[
		\E\sup_{h\in\H} (P-P_{n})h
		\le
		\frac{Ab}{n^{1/\alpha}}\left(
		1+ \frac{1}{(\alpha-2)} \left(\frac{c}{b} \right)^{\alpha/2}
		\right).
	\]	
\end{enumerate}
\end{lem}
\begin{proof} 
From the symmetrization principle, we get 
\begin{align*}
    \E\sup_{h\in\H} (P-P_{n})h
    \le 2R,
\end{align*}
where $R:= \E\sup_{h\in\H}R_n(h)$ for brevity.
Now it remains to bound $R$.
Consider the following cases of values $\alpha$.\\ \\
\noindent\textbf{\emph{Case of $\alpha<2$}}.
Here the square root of the entropy is still integrable around zero. Using \Cref{lem:chain} with $\eps=0$, we get
\begin{equation*}
	R \le \frac{12}{\sqrt n}\,\E \left[\int_{0}^{\sigma_n} \left(\frac{c}{u}\right)^{\alpha/2} \rmd u\right].
\end{equation*}
As in the proof of \Cref{lem:boundrad1}, 
using Jensen's inequality together with the symmetrization and contraction principles, we obtain
\begin{equation*}
	R \le \frac{12}{\sqrt n}\,\E \left[\int_{0}^{B} \left(\frac{c}{u}\right)^{\alpha/2} \rmd u\right],
\end{equation*}
where $B=\min\{b,\sigma+2\sqrt{bR}\}$. 
Integrating and using the fact that $B\ge\sigma$, yields
\begin{equation*}
	R \le \frac{24c^{\,\alpha/2}}{(2 - \alpha)\sqrt{n}} B^{1- \alpha/2} 
		\leq  \frac{24c^{\,\alpha/2}}{(2 - \alpha)\sqrt{n}}\sigma^{- \alpha/2} B 
		\leq  \frac{24c^{\,\alpha/2}}{(2 - \alpha)\sqrt{n}}\sigma^{- \alpha/2} (\sigma+2\sqrt{bR}).
\end{equation*}
Observing that $R\le x+y\sqrt R$ implies $R \le 2x + 2y^{2}$ for any $x,y>0$, we obtain 
\[
	R
	\leq
	\frac{48\sigma}{(2 - \alpha)\sqrt{n}}\left(\frac{c}{\sigma}\right)^{\alpha/2}
	+ \frac{4608b}{(2 - \alpha)^2\,n}\left(\frac{c}{\sigma}\right)^{\alpha}.
\]
\noindent\textbf{\emph{Case of $\alpha=2$}}.
Since the square root of the entropy is not anymore integrable around zero, we
apply \Cref{lem:chain} with $\eps=\sigma_n/n$ 
and obtain
\[
	R \leq \frac{4}{n} \E[\sigma_n]+ \frac{12c}{\sqrt n}\,\E \left[\int_{\sigma_n/n}^{\sigma_n} \frac{\rmd u}{u} \right] 
	\leq \frac{4b}{n} + \frac{12c\log n}{\sqrt n},
\]
where in the last inequality we have used the fact that, almost surely, $\sigma_n \leq b$. \\ \\
\noindent\textbf{\emph{Case of $\alpha>2$}}.
It follows from \Cref{lem:chain} that, for all $\eps > 0$, 
\[
	R
	\le 
	\E\left[4\eps+\frac{12}{\sqrt n}\int_{\eps}^{+\infty}\sqrt{\log \N\bigl(\H,\mn,u\bigr)}\rmd u\right].
\]
Taking $\eps = b/n^{1/\alpha}$ gives us 
\begin{align*}
	R 
	&\leq \frac{4b}{n^{1/\alpha}}+ \frac{12}{\sqrt n}\,\E \left[\int_{b/n^{1/\alpha}}^{+\infty} \left(\frac{c}{u}\right)^{\alpha/2} \rmd u\right]
	\leq \frac{4b}{n^{1/\alpha}} + \frac{24b}{(\alpha-2)n^{1/\alpha}} \left(\frac{c}{b} \right)^{\alpha/2}. 
\end{align*}
This finishes the proof.
\end{proof}

\subsection{Bounds on the expected supremum of empirical process using the bracketing entropy}\label{sec:bracketbounds}

In this section we derive explicit bounds for the expected supremum of empirical process 
under specific assumptions on the bracketing entropy of a class $\H$. 
We recall that the variance $\sigma_{\H}$ of the empirical process is defined as
\[
	\sigma_{\H}:=\sup_{h\in\H} \sqrt{Ph^2}.
\]
The following result is an analogue of Dudley's entropy bound, see \Cref{lem:chain}, but for the bracketing entropy.
This result follows from Lemma~7 in \citep{HWCS19} with the only difference that we present this result for $[-b,b]$-valued functions, 
not for $[-1,1]$-valued functions. However, a simple rescaling argument $\H \mapsto \H/b$ proves the claim. 
\begin{lem}
\label{lem:chainingbracket}
Let $\H$ be a class of $[-b,b]$-valued functions for some $b>0$. 
Then for any $\sigma > \sigma_{\H}$, any $\eps \in [0, \sigma/3]$, and any $n\ge1$, it holds that
\begin{align*}
	\E\sup_{h\in\H} \bigl|(P-P_{n})h\bigr|
	\le 
	A\left(
		\eps
		+\frac{1}{\sqrt{n}} \int_{\eps}^{\sigma}\sqrt{\log\Nb\bigl(\H, \m, u\bigr)} \rmd u 
		+ \frac{b}{n} \log \Nb\bigl(\H, \m, \sigma\bigr)
	\right),
\end{align*}
where $A>0$ is an absolute constant.
In particular, for $\sigma=2b$, we have 
\[
	\E\sup_{h\in\H} \bigl|(P-P_{n})h\bigr|
	\le 
	A\left(
		\eps
		+\frac{1}{\sqrt{n}} \int_{\eps}^{2b}\sqrt{\log\Nb\bigl(\H, \m, u\bigr)} \rmd u 
		+ \frac{b}{n}
	\right).
\]
\end{lem}

The next two lemmas apply \Cref{lem:chainingbracket} to derive explicit upper bounds on the
expected supremum of empirical process under specific assumptions on the bracketing metric entropy of $\H$.
The first case we study is whenever the class $\H$ has polynomial complexity.

\begin{lem}
\label{lem:boundrad3}
Let $\H$ be a class of $[-b,b]$-valued functions for some $b>0$. Suppose that, for all $0<u\le b$,  
\begin{equation*}
	\Nb\bigl(\H,\m,u\bigr)\le\left(\frac{c}{u}\right)^{\alpha},
\end{equation*}
for some constants 
$c>b\rme$ 
and $\alpha>0$. Then, for any  $\sigma \in [\sigma_{\H},b]$ and all $n\ge 1$,
\[
	\E\sup_{h\in\H} \bigl|(P-P_{n})h\bigr|
	\le
	A \left(\sqrt{\frac{\alpha\sigma^2}{n}\log\left(\frac{c}{\sigma}\right)} + \frac{\alpha b}{n}\log\left(\frac{c}{\sigma}\right)\right),
\]
where $A>0$ is an absolute constant. 
\end{lem}
\begin{proof}
Fix any $\sigma \in [\sigma_{\H},b]$.
It follows from \Cref{lem:chainingbracket} with $\eps = 0$ that 
\begin{align}
	\E\sup_{h\in\H} \bigl|(P-P_{n})h\bigr|
	\le 
	A\left(
		\frac{1}{\sqrt{n}} \int_{0}^{\sigma}\sqrt{\log\Nb\bigl(\H, \m, u\bigr)} \rmd u 
		+ \frac{b}{n} \log \Nb\bigl(\H, \m, \sigma\bigr)
	\right),
	\label{brentrop:e1}
\end{align}
where $A>0$ is an absolute constant.
Using the assumption on the bracketing number of $\H$, we get 
\begin{align*}
	\int_{0}^{\sigma}\sqrt{\log\Nb\bigl(\H, \m, u\bigr)} \rmd u 
	&\leq 
	\sqrt{\alpha} \int_{0}^{\sigma}\sqrt{\log \left(\frac{c}{u}\right)} \rmd u
	=
	\sqrt{\alpha} c \int_{c/\sigma}^{+\infty} u^{-2}\sqrt{\log \left(u\right)} \rmd u.
\end{align*}
Observe that, by assumptions, $c/\sigma\ge e$. Now, integration by parts yields
\[
	\sup_{t\ge e}\frac{\int_{t}^{+\infty} u^{-2}\sqrt{\log(u)}\,{\rm d}u}{t^{-1}\sqrt{\log(t)}}
	\leq 
	1 + \sup_{t\ge e}\frac{\frac12\int_{t}^{+\infty} u^{-2}(\log(u))^{-3/2}\,{\rm d}u}{t^{-1}\sqrt{\log(t)}}
	\leq
	\frac32.
\]
Therefore, we have
\begin{align*}
	\int_{0}^{\sigma}\sqrt{\log\Nb\bigl(\H, \m, u\bigr)} \rmd u 
	&\leq 
	\sqrt{\frac92 \alpha \sigma^2 \log \left(\frac{c}{\sigma}\right)}.
\end{align*}
Substituting this into \eqref{brentrop:e1}, for some new absolute constant $A>0$, we obtain 
\[
	\E\sup_{h\in\H} \bigl|(P-P_{n})h\bigr|
	\le
	A \left(\sqrt{\frac{\alpha\sigma^2}{n}\log\left(\frac{c}{\sigma}\right)} + \frac{\alpha b}{n}\log\left(\frac{c}{\sigma}\right)\right),
\]
which completes the proof.
\end{proof}

We turn to the case where the complexity of $\H$ is exponential. 
Similar to the metric entropy case, the square root of the bracketing entropy might not be integrable around zero anymore.

\begin{lem}
\label{lem:boundrad4}
Let $\H$ be a class of $[-b,b]$-valued functions for some $b>0$. Suppose that, for all $0<u\le 2b$,  
\begin{equation}
	\log\Nb\bigl(\H,\m,u\bigr)\le\left(\frac{c}{u}\right)^{\alpha},
\end{equation}
for some constants $c>0$ and $\alpha>0$. 
Then for some absolute constant $A>0$ and all $n\ge1$
the following holds.\\
\begin{enumerate}
	\item If $\alpha < 2$, then, for any  $\sigma \in [\sigma_{\H},b]$, 
	\[
		\E\sup_{h\in\H} \bigl|(P-P_{n})h\bigr|
		\le
		A\left(
		\frac{\sigma}{(2-\alpha)\sqrt{n}}\left(\frac{c}{\sigma}\right)^{\alpha/2}
		+ \frac{b}{n}\left(\frac{c}{\sigma}\right)^{\alpha}
		\right).
	\]
	\item  If $\alpha = 2$, then
	\[
		\E\sup_{h\in\H} \bigl|(P-P_{n})h\bigr|
		\le
		A\left(
		\frac{c \log n}{\sqrt{n}} + \frac{b}{n} 
		\right).
	\]	
	\item  If $\alpha > 2$, then
	\[
		\E\sup_{h\in\H} \bigl|(P-P_{n})h\bigr|
		\leq
		\frac{Ab}{n^{1/\alpha}}\left(
		1 + \frac{1}{(\alpha-2)} \left( \frac{c}{b}\right)^{\alpha/2}
		\right),
	\]	
\end{enumerate}
\end{lem}
\begin{proof}
Throughout the proof, let $A>0$ be a universal constant whose value may change from line to line. 
Consider the following cases. \\ \\
\noindent\textbf{\emph{Case of $\alpha<2$}}.
Fix any $\sigma \in [\sigma_{\H},b]$. 
It follows from \Cref{lem:chainingbracket} with $\eps=0$ that
\begin{align*}
	\E\sup_{h\in\H} \bigl|(P-P_{n})h\bigr|
	\le 
	A\left(
		\frac{1}{\sqrt{n}} \int_{0}^{\sigma}\sqrt{\log\Nb\bigl(\H, \m, u\bigr)} \rmd u 
		+ \frac{b}{n} \log \Nb\bigl(\H, \m, \sigma\bigr)
	\right).
\end{align*}
Using the assumption on the bracketing number of $\H$, we get 
\begin{align*}
	\int_{0}^{\sigma}\sqrt{\log\Nb\bigl(\H, \m, u\bigr)} \rmd u 
	&\leq 
	\int_{0}^{\sigma} \left(\frac{c}{u}\right)^{\alpha/2} \rmd u
	=
	\frac{2c^{\alpha/2}\,\sigma^{1-\alpha/2}}{2-\alpha}.
\end{align*}
Therefore, we conclude 
\begin{align*}
	\E\sup_{h\in\H} \bigl|(P-P_{n})h\bigr|
	\leq 
	A\left(
		\frac{\sigma}{(2-\alpha)\sqrt{n}}\left(\frac{c}{\sigma}\right)^{\alpha/2}
		+ \frac{b}{n}\left(\frac{c}{\sigma}\right)^{\alpha}
	\right).
\end{align*}
\noindent\textbf{\emph{Case of $\alpha=2$}}.
Since the square root of the bracketing entropy is not anymore integrable around zero, we
apply \Cref{lem:chainingbracket} with $\sigma=2b$ and $\eps=b/n$ and obtain for $n\ge3/2$
\begin{align*}
	\E\sup_{h\in\H} \bigl|(P-P_{n})h\bigr|
	&\leq 
	A\left(
		\frac{b}{n} + \frac{1}{\sqrt{n}} \int_{b/n}^{2b} \frac{c}{u} \rmd u + \frac{b}{n}
	\right)
	\leq 
	A\left(
		\frac{b}{n} + \frac{c \log n}{\sqrt{n}} 
	\right).
\end{align*}
It is clear that the same bound holds for $n=1$. \\ \\
\noindent\textbf{\emph{Case of $\alpha>2$}}.
Using  \Cref{lem:chainingbracket} with $\sigma=2b$ and $\eps=b/n^{1/\alpha}$, we obtain for $n\ge3/2$
\begin{align*}
	\E\sup_{h\in\H} \bigl|(P-P_{n})h\bigr|
	&\leq 
	A\left(
		\frac{b}{n^{1/\alpha}} + \frac{1}{\sqrt{n}} \int_{b/n^{1/\alpha}}^{+\infty} \left(\frac{c}{u}\right)^{\alpha/2} \rmd u + \frac{b}{n}
	\right)
    \\
	&\leq 
	A\left(
		\frac{b}{n^{1/\alpha}} + \frac{2 b}{(\alpha-2)n^{1/\alpha}} \left( \frac{c}{b}\right)^{\alpha/2}
	\right).
\end{align*}
The same bound obviously holds for $n=1$. This finishes the proof of the lemma.
\end{proof}
\subsection{Proof of \texorpdfstring{\Cref{lem:fixedpoint}}{Lg}}\label{sec:lemfixedpoint}
As it was mentioned, this result follows from a combination of arguments presented in \citep{Kol11}. 
We devide the proof into two steps.
\\ \\
\noindent\textbf{\emph{Step 1}}.
Let $\delta_{j}$, $j\ge 0$, be a strictly decreasing sequence of positive numbers with $\delta_{0}=\bar \delta$ and let $t_{j}$, $j\ge 0$, be an arbitrary sequence of positive numbers. For all $\delta\ge0$, denote
\[
	\bar\beta(\delta)=\sum_{j=0}^{+\infty}\beta\left(\delta_{j},t_{j}\right)\mathbf 1\left\{\delta_{j+1}<\delta\le \delta_{j}\right\},
\]
and set 
\[
	\delta^*=\sup\left\{\delta\ge0:\delta\le \bar \beta(\delta)\right\}.
\]
The goal of this first step is to prove that, for all $\delta\ge \delta^*$, 
\[
	\P(\hat \delta\ge \delta)\le \sum_{j:\,\delta_{j}\ge \delta}e^{-t_{j}}.
\]
Fix any $\delta>\delta^*$. For all $j\ge 0$, define  $E_{j}=\left\{\phi(\delta_{j})< \bar \beta(\delta_{j})\right\}$ and set
\[
	E=\bigcap _{j:\,\delta_{j}\ge \delta}E_{j}.
\]
Combining property \eqref{lem:fixedpoint:e1} with the definition of $\bar\beta$ 
(in particular that $\bar\beta(\delta_j)=\beta(\delta_j,t_j)$ by construction) yields
\[
	\P(E)\ge 1-\sum_{j:\,\delta_{j}\ge \delta}e^{-t_{j}}.
\] 
On the event $E$ and for all $\delta'\ge \delta$, we have $\phi(\delta')\le \bar \beta(\delta')$, 
by the monotonicity of $\phi$ and by the definition of $\bar \beta$. Thus, on the event $\{\hat \delta\ge \delta\}\cap E$ we obtain
\[
	\hat \delta\le \phi(\hat \delta)\le \bar \beta(\hat \delta),
\]
which implies that $\delta\le \hat \delta\le \delta^*$. Since this contradicts $\delta>\delta^*$, we deduce that $\{\hat \delta\ge \delta\}\subset E^c$ and hence
\[
	\P(\hat \delta\ge\delta)< \sum_{j:\,\delta_{j}\ge \delta}e^{-t_{j}}.
\]
By continuity, this also holds for $\delta=\delta^*$.
\\ \\
\noindent\textbf{\emph{Step 2}}. 
Fix $t>0$ and set $\tau = \beta(t)$ with $\beta(t)$ defined in \eqref{lem:fixedpoint:e2}
(there is no loss of generality in assuming that $\beta(t)$ is finite).
Set for $j\ge0$, 
\[
    \delta_j = \frac{\bar\delta}{2^{j}}
    \quad\text{and}\quad
    t_j = \frac{t\bar\delta }{\tau 2^{j}}.
\]
Now for this specific choice of $\delta_j$ and $t_j$ observe that, for any $\delta\ge\tau$,
\begin{align*}
	\frac{\bar \beta(\delta)}{\delta} 
	&= \sum_{j=0}^{+\infty}\frac{\bar\delta}{\delta 2^{j}}\left(\frac{\bar\delta}{2^{j}}\right)^{-1}\beta\left(\frac{\bar \delta}{2^{j}},\frac{t\bar \delta}{\tau2^{j}}\right)\mathbf 1\left\{\frac{\bar \delta}{2^{j+1}}<\delta\le\frac{\bar \delta}{2^{j}}\right\}\\
	&\le 2 \sum_{j=0}^{+\infty}\left(\frac{\bar\delta}{2^{j}}\right)^{-1}\beta\left(\frac{\bar \delta}{2^{j}},\frac{t\bar \delta}{\tau2^{j}}\right)\mathbf 1\left\{\frac{\bar \delta}{2^{j+1}}<\delta\le\frac{\bar \delta}{2^{j}}\right\}\\
	&\le 2 \sup_{\delta\ge \tau}\frac{\beta\left(\delta,\tfrac{t\delta}{\tau}\right)}{\delta}\\
	&\le1,
\end{align*}
where the last inequality follows from the definition of $\beta(t)$.
Hence
\[
	\delta\ge \delta^{\star}=\sup\left\{\delta\ge 0:\, \delta\le \bar \beta(\delta)\right\}.
\]
Then, according to the first step above, 
\begin{equation}
\label{ftr:e2}
	\P(\hat \delta\ge \delta)\le \sum_{j:\,\frac{\bar \delta}{2^j}\ge \delta}e^{-\frac{t\bar \delta}{\tau 2^j}}.
\end{equation}
The sum on the right-hand side of \eqref{ftr:e2} may be bounded as follows. Let 
\[
	j^*=\max\left\{j\ge 0:\frac{\bar \delta}{2^j}\ge \tau\right\}.
\]
Then
\begin{equation*}
	\sum_{j:\frac{\bar \delta}{2^j}\ge \delta}e^{-\frac{t\bar \delta}{\tau 2^j}} 
	=\sum_{j=0}^{j^*}e^{-\frac{t\bar \delta}{\tau 2^j}}\le \sum_{j=0}^{j^{*}}e^{-t2^{(j^{*}-j)}}\le\sum_{j=0}^{+\infty}e^{-t 2^{j}},
\end{equation*}
and
\begin{equation*}
	\sum_{j=0}^{+\infty}e^{-t 2^{j}}\le e^{-t}+\sum_{j=1}^{+\infty}(2^j-2^{j-1})e^{-t 2^{j}}\le e^{-t}+\int_{1}^{+\infty}e^{-tu}{\rm d}u= 2e^{-t}.
\end{equation*}
Finally, we have proved that, for all $t\ge 0$ and for all $\delta>\tau=\beta(t)$ we have 
\[
	\P\bigl(\hat \delta\ge \delta\bigr)\le 2e^{-t}.
\]
By continuity, this result also holds for $\delta=\beta(t)$.
\end{appendix}

\end{document}